\documentclass[a4paper,10pt]{article}

\usepackage{latexsym}
\usepackage{amssymb}
\usepackage{theorem}
\usepackage{amsmath}
\usepackage{amscd}
\usepackage{color}
\usepackage{graphicx}

\newtheorem{theorem}{Theorem}[section]
\newtheorem{corollary}[theorem]{Corollary}
\newtheorem{lemma}[theorem]{Lemma}
\newtheorem{example}[theorem]{Example}
\newtheorem{proposition}[theorem]{Proposition}
\newtheorem{remark}[theorem]{Remark}
\newtheorem{definition}[theorem]{Definition}

\newcommand{\demo}{\par\noindent{\it Proof. \/}\ }
\newcommand{\enD}{\hfill $\Box$\vspace{3truemm} \par}
\newcommand{\R}{\mathbb{R}}

\newcommand{\bX}{\mbox{\boldmath $X$}}
\newcommand{\ba}{\mbox{\boldmath $a$}}
\newcommand{\bx}{\mbox{\boldmath $x$}}
\newcommand{\by}{\mbox{\boldmath $y$}}
\newcommand{\bz}{\mbox{\boldmath $z$}}
\newcommand{\bt}{\mbox{\boldmath $t$}}
\newcommand{\be}{\mbox{\boldmath $e$}}
\newcommand{\bv}{\mbox{\boldmath $v$}}
\newcommand{\bw}{\mbox{\boldmath $w$}}
\newcommand{\bm}{\mbox{\boldmath $m$}}
\newcommand{\bn}{\mbox{\boldmath $n$}}
\newcommand{\bs}{\mbox{\boldmath $s$}}
\newcommand{\bF}{\mbox{\boldmath $F$}}
\newcommand{\bV}{\mbox{\boldmath $V$}}
\newcommand{\bN}{\mbox{\boldmath $N$}}

\begin{document}

\title{Focal surfaces of lightcone framed surfaces in the Lorentz-Minkowski $3$-space}

\author{Chang Xu and Liang Chen}
\date{}

\maketitle

\begin{abstract}

In this paper, we consider the differential geometry properties of focal surfaces of lightcone framed surfaces in Lorentz-Minkowski $3$-space.
In general, a mixed type surface is a connected regular surface with non-empty spacelike and timelike point sets.
In order to investigate the mixed type surface with singular points, we introduce the lightcone framed surface.
First, we give the lightcone frame, by which we define the lightcone framed surface.
Next, we consider the differential geometry properties of lightcone framed surfaces by using the lightcone frame.
At last, we give the definition of focal surfaces of lightcone framed surfaces and investigate the differential geometry properties of focal surfaces.

\end{abstract}

\section{Introduction}

If a point of the surface moves between spacelike and timelike regions, then there is a lightlike point. 
As we all know, a mixed type surface is a connected regular surface in a Lorentzian $3$-manifold with non-empty spacelike and timelike
point sets.
The induced metric of a mixed type surface is a signature-changing metric, and their lightlike points can be regarded as singular points of this metric.
So lightlike points play key roles in studying a mixed type surface.
In this paper, we consider lightlike points of lightcone framed surfaces in the way similar to the case of singular points of fronts.
First we define the singular mixed type surface, so-called lightcone framed surface, with not only lightlike points but also singular points.
In order to investigate the lightcone framed surface, we construct the lightcone frame.
The idea is that the framed surface in the Euclidean $3$-space \cite{Fukunaga1} and the lightcone framed curve in Lorentz-Minkowski $3$-space \cite{Liang Chen1}.
Then we consider the differential geometric properties of the lightcone framed surface by using the lightcone frame and classify the lightlike points of the lightcone framed surface by using the criterion of singular points of fronts.

The focal surface is a classical object in differential geometry, which can be characterized in several way: loci of centers of principal curvature spheres of initial surface, singular value sets of the normal congruence and the bifurcation set of the family of distance squared functions for instance, for details see \cite{Tari1, Teramoto1, Teramoto2, Teramoto3, Teramoto4}.
It is known that focal surfaces of regular surfaces have singular points in general.
Principal curvature and principal direction of surfaces are important concepts in differential geometry of surfaces and they are closely related to focal surfaces.
So studying focal surfaces can provide us with new geometric properties of surfaces.
However, focal surfaces of mixed type surfaces may have not only singular points but also lightlike points.
F.Tari considered focal surfaces of regular mixed type surfaces in Lorentz-Minkowski $3$-space by using the bifurcation set of the family of distance squared functions (see \cite{Tari1}).
In this paper, we consider focal surfaces of lightcone framed surfaces and give the differential geometry properties of focal surfaces.

The organization of this paper is as follows.
In section $2$, we prepare some basic notions of Lorentz-Minkowski $3$-space (cf. \cite{Liang Chen1} \cite{T.Liu1}).
Morever we review on lightcone frame of mixed type curves in Lorentz-Minkowski $3$-space (cf. \cite{Liang Chen1}) and framed surfaces in the Euclidean $3$-space (cf. \cite{Fukunaga1}).
In section $3$, we define the lightcone framed surface in Lorentz-Minkowski $3$-space.
Next we discuss the Gaussian curvature and the Mean curvature of the lightcone framed surface not only in non-lightlike part but also with lightlike points.
Futhermore, similar to singular points of wave fronts in the Euclidean $3$-space, we investigate the properties of lightlike points of lightcone framed surfaces in Lorentz-Minkowski $3$-space.
In section $4$, we give the principal curvature and the principal direction of the lightcone framed surface and the definition of the focal surface of the lightcone framed surface in Lorentz-Minkowski $3$-space by using the bifurcation set of the family of distance-squared functions of the surface.
Morever, similar to the case of the lightcone framed surface, we discuss the Gaussian curvature and the Mean curvature of the focal surface of the lightcone framed of not only surface on non-lightlike points but also surface with lightlike points.
At last, we consider the differential geometry properties of the Gaussian curvature and the Mean curvature of the focal surface of the lightcone framed surface.

We shall assume throughout the whole paper that all maps and manifolds are differentiable of class $C^\infty$ and arguments of functions of surfaces are all $(u,v)$ without unless stated.

\section{Preliminaries}

\subsection{Lorentz-Minkowski $3$-space}

Let $\R^3=\{ (x_1,x_2,x_3)|x_i\in \R, i=1,2,3 \}$ be a $3$-dimensional vector space. 
For any vectors $\bx=(x_1,x_2,x_3),\by=(y_1,y_2,y_3)$ in $\R^3$, the {\it pseudo inner product} of $\bx$ and $\by$ is defined to be 
$\langle \bx,\by \rangle=-x_1y_1+x_2y_2+x_3y_3$. 
We call $(\R^3,\langle , \rangle)$ the {\it Lorentz-Minkowski 3-space} and write $\R_1^3$ instead of $(\R^3,\langle ,\rangle)$.

We say that non-zero vector $\bx$ in $\R_1^3$ is {\it spacelike, lightlike} or {\it timelike} if $\langle \bx,\bx \rangle>0$, $ \langle \bx,\bx \rangle=0$ or $\langle \bx,\bx \rangle<0$ repectively. 
The {\it norm} of the vector $\bx$ is defined by $ \| \bx \| =\sqrt{| \langle \bx,\bx\rangle|}$.

For any $\bx=(x_1,x_2,x_3)$, $\by=(y_1,y_2,y_3)\in \R_1^3$, we define a vector $\bx\wedge\by$ by 
\begin{eqnarray*}
\bx\wedge\by=\begin{vmatrix}-\be_1
&\be_2  &\be_3 \\x_1
&x_2  &x_3 \\y_1
&y_2  &y_3
\end{vmatrix},
\end{eqnarray*}
where $\{\be_1, \be_2, \be_3\}$ is the canonical basis of $\R_1^3$. 

For any $\bz \in \R_1^3$, we can easily check that 
\begin{eqnarray*}
\langle \bz,\bx\wedge\by \rangle=det(\bz,\bx,\by),
\end{eqnarray*}
so that $\bx \wedge\by$ is pseudo-orthogonal to both $\bx$ and $\by$, and
\begin{eqnarray*}
&&(\bx\wedge\by)\wedge\bz=-\langle \bx,\bz \rangle\by+\langle \by,\bz \rangle\bx,\\
&&\bx\wedge(\by\wedge\bz)=\langle \bx,\by \rangle\bz-\langle \bx,\bz \rangle\by.
\end{eqnarray*}
Morever for any $\ba\in\R_1^3$, we have
\begin{eqnarray*}
(\bx\wedge\by,\bz\wedge\ba)=\langle \bx,\ba \rangle \langle\by,\bz\rangle-\langle \bx,\bz \rangle \langle\by,\ba\rangle.
\end{eqnarray*}

If $\bx$ is a timelike vector, $\by$ is a spacelike vector and $\bx\wedge\by=\bz$, then by a straightforward calculation, we have 
\begin{eqnarray*}
\bz\wedge\bx=\by, \by\wedge\bz=-\bx.
\end{eqnarray*}
If $\bx$ is a spacelike vector, $\by$ is a timelike vector and $\bx\wedge\by=\bz$, then by a straightforward calculation, we have 
\begin{eqnarray*}
\bz\wedge\bx=-\by, \by\wedge\bz=\bx.
\end{eqnarray*}
If both $\bx$ and $\by$ are spacelike vectors and $\bx\wedge\by=\bz$, then by a straightforward calculation, we have 
\begin{eqnarray*}
 \bz\wedge\bx=-\by, \by\wedge\bz=-\bx.
\end{eqnarray*}

For a vector $\bv \in \R_1^3$ and a real number c, we define the {\it plane} with the pseudo-normal $\bv$ by 
\begin{eqnarray*}
 P(\bv,c)=\left \{ \bx\in \mathbb{R}_1^3|\left \langle  \bx,\bv\right \rangle=c \right \}.
\end{eqnarray*}
We call $P(\bv,c)$ a {\it timelike plane}, {\it lightlike plane} or {\it spacelike plane} if $\bv$ is spacelike, lightlike or timelike, respectively.

We define {\it hyperbolic $2$-space} by
\begin{eqnarray*}
H^2(-1)=\{\bx\in \R_1^3\ |\ \langle\bx,\bx\rangle =-1\},
\end{eqnarray*}
{\it de Sitter $2$-space} by
\begin{eqnarray*}
S^2_1=\{\bx\in \R_1^3\ |\ \langle\bx,\bx\rangle =1\},
\end{eqnarray*}
{\it (open) lightcone} at the origin by
\begin{eqnarray*}
LC^{*}=\{\bx \in \R_1^3\setminus \{0\} \ |\ \langle \bx, \bx\rangle=0\}.
\end{eqnarray*}

In this paper, we consider a double Legendrian fibration :
\begin{itemize}
\item[(1)]
$\Delta_4 =\{(\bv,\bw) \in LC^* \times LC^* \ | \ \langle \bv,\bw \rangle=-2 \},$
\par
$\pi_{41}:\Delta_4 \to LC^*, \pi_{41}(\bv,\bw)=\bv, \ \pi_{42}:\Delta_4 \to LC^*, \pi_{42}(\bv,\bw)=\bw$, 
\par
$\theta_{41}=\langle d\bv, \bw \rangle|_{\Delta_4}, \theta_{42}=\langle \bv,d\bw \rangle|_{\Delta_4}$.
\end{itemize}

\subsection{Lightcone frame of mixed type curves in the Lorentz-Minkowski $3$-space}

Let $I$ be an interval of $\R$ and $\gamma:I\to\R_1^3$ be a smooth curve.
We say that $\gamma$ is {\it spacelike} (respectively, {\it lightlike, timelike}) if $\dot{\gamma}(t)=(d\gamma/dt)(t)$ is a spacelike (respectively, lightlike, timelike) vector for any $t \in I$. 
We call $\gamma$ a {\it non-lightlike} curve if $\gamma$ is spacelike or timelike.
Moreover, a point $t$ (or, $\gamma(t)$) is called a {\it spacelike} (respectively, {\it lightlike}, {\it timelike}) point if $\dot{\gamma}(t)$ is a spacelike (respectively, lightlike, timelike) vector. 

Let $(\gamma,\ell^+,\ell^-):I \to \R_1^3\times\Delta_4$ be a smooth mapping. 
We say that $(\gamma,\ell^+,\ell^-)$ is a {\it lightcone framed curve} if there exists smooth functions $\alpha,\beta:I \to \R$ such that $\dot{\gamma}(t)=\alpha(t)\ell^+(t)+\beta(t)\ell^-(t)$ for all $t \in I$. 
We also say that $\gamma$ is a {\it lightcone framed base curve} if there exists a smooth mapping $(\ell^+,\ell^-):I\to \R_1^3\times\Delta_4$ such that $(\gamma,\ell^+,\ell^-):I\to \R_1^3\times\Delta_4$ is a lightcone framed curve.

We say that $\gamma$ is {\it spacelike, lightlike} or {\it timelike} at $t$ if $\alpha(t)\beta(t)<0, \alpha(t)\beta(t)=0$ or $\alpha(t)\beta(t)>0$, respectively. 
Note that if $t$ is a lightlike point of $\gamma$, then $\alpha(t)\beta(t)=0$ with $\alpha(t)\ne0$ or $\beta(t)\ne0$. 
Moreover, $t$ is a singular point of $\gamma$ if and only if $\alpha(t)=\beta(t)=0$.

Let $(\gamma,\ell^+,\ell^-):I \to \R_1^3\times\Delta_4$ be a lightcone framed curve and $\bn(t)=-(1/2)\ell^+(t)\wedge\ell^-(t)$. 
We say that $\{\ell^+(t),\ell^-(t),\bn(t)\}$ is a {\it lightcone frame} of $\gamma(t)$ and
\begin{eqnarray*}
\left(
\begin{array}{c}
\dot{\ell^{+}}(t)\\
\dot{\ell^{-}}(t)\\
\dot{\bn}(t)
\end{array}
\right)
&=&
\left(
\begin{array}{ccc}
\kappa_1(t)&0&2\kappa_3(t) \\
0&-\kappa_1(t)&2\kappa_2(t)\\
\kappa_2(t)&\kappa_3(t)&0
\end{array}
\right)
\left(
\begin{array}{c}
\ell^+(t)\\
\ell^-(t)\\
\bn(t)
\end{array}
\right), \\
\dot{\gamma}(t)&=&\alpha(t)\ell^+(t)+\beta(t)\ell^-(t),
\end{eqnarray*}
where 
\begin{eqnarray*}
&&\kappa_1(t)=-\frac{1}{2} \langle \dot{\ell^+}(t),\ell^-(t) \rangle, \ 
\kappa_2(t)=-\frac{1}{2} \langle \dot{\bn}(t),\ell^-(t) \rangle, \ 
\kappa_3(t)=-\frac{1}{2} \langle \dot{\bn}(t),\ell^+(t) \rangle, \\
&&\alpha(t)=-\frac{1}{2} \langle \dot{\gamma}(t), \ell^-(t) \rangle, \ 
\beta(t)=-\frac{1}{2} \langle \dot{\gamma}(t), \ell^+(t) \rangle.
\end{eqnarray*}

We denote $\overline\ell^+(t)=(1/c(t))\ell^+(t)$, $\overline\ell^-(t)=c(t)\ell^-(t)$, where $c:I \to \R$ is a non-zero function. 
Then we get $\overline\ell^+(t)\wedge\overline\ell^-(t)=\ell^+(t)\wedge\ell^-(t)=-2\bn(t)$. 
So $\{\overline\ell^+(t),\overline\ell^-(t),\bn(t) \}$ is also a lightcone frame of $\gamma(t)$. 
By a direct calculation, we have
\begin{eqnarray*}
\left(
\begin{array}{c}
\dot{\overline{\ell}^{+}}(t)\\
\dot{\overline{\ell}^{-}}(t)\\
\dot{\bn}(t)
\end{array}
\right)
&=&
\left(
\begin{array}{ccc}
\overline{\kappa}_1(t)&0&2\overline{\kappa}_3(t) \\
0&-\overline{\kappa}_1(t)&2\overline{\kappa}_2(t)\\
\overline{\kappa}_2(t)&\overline{\kappa}_3(t)&0
\end{array}
\right)
\left(
\begin{array}{c}
\overline{\ell}^+(t)\\
\overline{\ell}^-(t)\\
\bn(t)
\end{array}
\right), \\
\dot{\gamma}(t)&=&\overline{\alpha}(t) \overline{\ell}^+(t)+\overline{\beta}(t) \overline{\ell}^-(t),
\end{eqnarray*}
where 
\begin{eqnarray*}
\overline{\kappa}_1(t) =\frac{c(t)\kappa_1(t)-\dot{c}(t)}{c(t)}, \ \overline{\kappa}_2(t)=c(t)\kappa_2(t), \ 
\overline{\kappa}_3(t)=\frac{\kappa_3(t)}{c(t)}, \ 
\overline{\alpha}(t)=c(t)\alpha(t), \ 
\overline{\beta}(t)=\frac{\beta(t)}{c(t)}.
\end{eqnarray*}
If we take $\dot{c}(t)=c(t)\kappa_1(t)$, that is, $c(t)=Ae^{\int \kappa_1(t) dt}$, where $A$ is a constant, then we get $\overline{\kappa}_1(t)=0$.
We say that the lightcone frame $\{ \overline{\ell}^+(t),\overline{\ell}^-(t),\bn(t) \}$ with $\overline{\kappa}_1(t)=0$ is an {\it adapted frame}.  
We rewire  $(\gamma,\ell^+,\ell^-):I \to \R_1^3\times\Delta_4$ with curvature $(m,n,\alpha,\beta)$, that is,
\begin{eqnarray*}
\left(
\begin{array}{c}
\dot{\ell^{+}}(t)\\
\dot{\ell^{-}}(t)\\
\dot{\bn}(t)
\end{array}
\right)
&=&
\left(
\begin{array}{ccc}
0&0&2n(t) \\
0&0&2m(t)\\
m(t)&n(t)&0
\end{array}
\right)
\left(
\begin{array}{c}
\ell^+(t)\\
\ell^-(t)\\
\bn(t)
\end{array}
\right), \\
\dot{\gamma}(t)&=&\alpha(t)\ell^+(t)+\beta(t)\ell^-(t).
\end{eqnarray*}
We call $(m,n,\alpha,\beta)$ a {\it lightcone curvature} of the lightcone framed curve $(\gamma,\ell^+,\ell^-):I\to \R_1^3\times\Delta_4$.

\subsection{Framed surfaces in the Euclidean $3$-space}

Let $\R^3$ be a $3$-dimensional Euclidean space equipped with the inner product $\bx\cdot\by=x_1y_1+x_2y_2+x_3y_3$, where $\bx=(x_1,x_2,x_3),\by=(y_1,y_2,y_3),\bx,\by \in \R$.
The norm of $\bx$ is given by $|\bx|=\sqrt{\bx\cdot\bx}$ and the vector product is given by
\begin{eqnarray*}
\bx \times \by=\begin{vmatrix}
  \be_1&  \be_2& \be_3\\
  x_1&  x_2& x_3\\
  y_1&  y_2& y_3
\end{vmatrix},
\end{eqnarray*}
where $\{ \be_1,\be_2,\be_3 \}$ is the canonical basis on $\R^3$. 
Let $U$ be a connected domain on $\R^2$ and $S^2=\{ \bx\in\R^3| |\bx|=1\}$ be a unit sphere in $\R^3$.
We denote a $3$-dimensional manifold by $\Delta=\{ (\bx,\by)\in S^2\times S^2| \bx\cdot\by=0 \}$.
Then we define $(\bx,\bn,\bs):U\to \R^3\times \Delta$ be a {\it framed surface} if $\bx_u\cdot\bn=\bx_v\cdot\bn=0$ for all $(u,v)\in U$, where $\bx_u=(\partial/\partial_u),\bx_v=(\partial/\partial_v)$.
We say that $\bx:U\to \R^3$ is a {\it framed base surface} if there exists $(\bn,\bs):U \to \Delta$ such that $(\bx,\bn,\bs)$ is a framed surface.

We denote $\bt=\bn\times\bs$, then $\{\bn,\bs,\bt\}$ is a moving frame along $\bx$.
By a direct calculation, we have the following equations:
\begin{eqnarray}\label{formula1}
\begin{pmatrix}
 \bx_u\\\bx_v
\end{pmatrix}=\begin{pmatrix}
  a_1& b_1\\
  a_2& b_2
\end{pmatrix}\begin{pmatrix}
 \bs\\\bt
\end{pmatrix},
\end{eqnarray}

\begin{eqnarray}\label{formula2}
\begin{pmatrix}
 \bn_u\\
 \bs_u\\
 \bt_u
\end{pmatrix}=\begin{pmatrix}
  0&  e_1&f_1 \\
  -e_1&  0&g_1 \\
  -f_1&-g_1&0
\end{pmatrix}\begin{pmatrix}
 \bn\\
 \bs\\
 \bt
\end{pmatrix},
\begin{pmatrix}
 \bn_v\\
 \bs_v\\
 \bt_v
\end{pmatrix}=\begin{pmatrix}
  0&  e_2&f_2 \\
  -e_2&  0&g_2 \\
  -f_2&-g_2&0
\end{pmatrix}\begin{pmatrix}
 \bn\\
 \bs\\
 \bt
\end{pmatrix},
\end{eqnarray}
where $a_i,b_i,e_i,f_i,g_i:U\to\R,i=1,2$ are smooth functions and we call them {\it basic invariants} of the framed surface.
We denote the matrices in \eqref{formula1} and \eqref{formula2} by $\mathcal{G},\mathcal{F}_1,\mathcal{F}_2$, respectively.
We call the matrices $(\mathcal{G},\mathcal{F}_1,\mathcal{F}_2)$ {\it basic invariants} of the framed surface $(\bx,\bn,\bs)$.
Note that $(u,v)$ is a singular point of $\bx$ if and only if $det\mathcal{G}=0$.

By the integrability conditions $\bx_{uv}=\bx_{vu},\mathcal{F}_{2u}-\mathcal{F}_{1v}=\mathcal{F}_1\mathcal{F}_2-\mathcal{F}_2\mathcal{F}_1$, we have
\begin{eqnarray*}
\left\{\begin{matrix}
 a_{1v}-b_1g_2=a_{2u}-b_2g_1,\\
 b_{1v}-a_2g_1=b_{2u}-a_1g_2,\\
 a_1e_2+b_1f_2=a_2e_1+b_2f_1,
\end{matrix}\right. \\
\left\{\begin{matrix}
 e_{1v}-f_1g_2=e_{2u}-f_2g_1,\\
 f_{1v}-e_2g_1=f_{2u}-e_1g_2,\\
 g_{1v}-e_1f_2=g_{2u}-e_2f_1.
\end{matrix}\right.
\end{eqnarray*}

\section{Lightcone framed surfaces}

\subsection{Lightcone framed surfaces}

In order to investigate mixed type of surfaces in the Lorentz-Minkowski space, we introduce the lightcone framed surface. 

\begin{definition}[\cite{M.Li1}]
Let $(\bX,\bv,\bw):U \to \R^3_1 \times \Delta_4$ be a smooth mapping.
We say that $(\bX,\bv,\bw):U \to \R^3_1 \times \Delta_4$ is a {\it lightcone framed surface} if there exist smooth functions: $\alpha,\beta: U\to \R$ such that 
\begin{eqnarray*}
\bX_u\wedge\bX_v=\alpha\bv+\beta\bw
\end{eqnarray*}
for all $(u,v)\in U$.
We also say that $\bX:U \to \R^3_1$ is a {\it lightcone framed base surface} if there exists a smooth mapping $(\bv,\bw):U \to \Delta_4$ such that $(\bX,\bv,\bw):U \to \R^3_1 \times \Delta_4$ is a lightcone framed surface.
\end{definition}

We denote $\bm=-(1/2)\bv \wedge \bw$.
It is obvious that $\langle\bm,\bm\rangle=1$.
We call $\{\bv,\bw,\bm\}$ the {\it lightcone frame} of $\Bbb R_1^3$ along $\bX$.
By a direct calculation, we have the following differential equations:

\begin{eqnarray}\label{curvature-one}
	\begin{pmatrix}
		\bX_u\\\bX_v
	\end{pmatrix}=\begin{pmatrix}
		a_1&  b_1& c_1\\
		a_2&  b_2& c_2
	\end{pmatrix}\begin{pmatrix}
		\bv\\
		\bw\\
		\bm
	\end{pmatrix}
\end{eqnarray}

\begin{eqnarray}\label{curvature-two}
	\begin{pmatrix}
		\bv_u\\\bw_u\\\bm_u
	\end{pmatrix}=\begin{pmatrix}
		e_1&  0& 2g_1\\
		0& -e_1& 2f_1\\
		f_1& g_1& 0
	\end{pmatrix}\begin{pmatrix}
		\bv\\
		\bw\\
		\bm
	\end{pmatrix},\begin{pmatrix}
		\bv_v\\\bw_v\\\bm_v
	\end{pmatrix}=\begin{pmatrix}
		e_2&  0& 2g_2\\
		0& -e_2& 2f_2\\
		f_2& g_2& 0
	\end{pmatrix}\begin{pmatrix}
		\bv\\
		\bw\\
		\bm
	\end{pmatrix}
\end{eqnarray}

where 
\begin{eqnarray*}
	&&a_1=-\frac{1}{2}\langle \bX_u,\bw \rangle,b_1=-\frac{1}{2}\langle \bX_u,\bv \rangle,c_1=\langle \bX_u,\bm \rangle,\\
	&&a_2=-\frac{1}{2}\langle \bX_v,\bw \rangle,b_2=-\frac{1}{2}\langle \bX_v,\bv \rangle,c_2=\langle \bX_v,\bm \rangle,\\
	&&e_1=\frac{1}{2}\langle \bv,\bw_u \rangle,f_1=\frac{1}{2}\langle \bw_u,\bm \rangle,g_1=\frac{1}{2}\langle \bv_u,\bm \rangle,\\
	&&e_2=\frac{1}{2}\langle \bv,\bw_v \rangle,f_2=\frac{1}{2}\langle \bw_v,\bm \rangle,g_2=\frac{1}{2}\langle \bv_v,\bm \rangle.
\end{eqnarray*}
Note that $a_i,b_i,c_i,e_i,f_i,g_i:U \to \R,i=1,2$ are smooth functions and we call them {\it basic invariants} of the lightcone framed surface $(\bX,\bv,\bw)$.
We denote the matrices given in \eqref{curvature-one} and \eqref{curvature-two} by $\mathcal G,\mathcal F_1,\mathcal F_2$, respectively.
We also call the matrices $(\mathcal G,\mathcal F_1,\mathcal F_2)$ {\it basic invariants} of the lightcone framed surface $(\bX,\bv,\bw)$.
Since the integrability conditions $\bX_{uv}=\bX_{vu}$ and $\mathcal F_{2u}-\mathcal F_{1v}=\mathcal F_1\mathcal F_2-\mathcal F_2\mathcal F_1$, basic invariants satisfy the following conditions:

\begin{eqnarray}
	\left\{\begin{matrix}
		a_{1v}=-a_1e_2+c_2f_1-c_1f_2\\
		b_{1v}=b_1e_2+c_2g_1-c_1g_2\\
		c_{1v}-c_{2u}=-2(a_1g_2+b_1f_2)
	\end{matrix}\right.\\
	\left\{\begin{matrix}
		e_{1v}-e_{2u}=2(f_1g_2-f_2g_1)\\
		f_{1v}-f_{2u}=e_1f_2-e_2f_1\\
		g_{1v}-g_{2u}=e_2g_1-e_1g_2
	\end{matrix}\right.
\end{eqnarray}

\begin{lemma}
	Under the settings above, we have $\bv\wedge\bm=-\bv, \bw\wedge\bm=\bw$.
\end{lemma}
\demo
Since $\bm=-(1/2)\bv\wedge\bw$, we have
\begin{eqnarray*}
	&&\bv\wedge\bm=\bv\wedge(-(1/2)\bv\wedge\bw)=-(1/2)(\langle \bv,\bv \rangle\bw-\langle \bv,\bw \rangle\bv)=-\bv,\\
	&&\bw\wedge\bm=\bw\wedge(-(1/2)\bv\wedge\bw)=-(1/2)(\langle \bv,\bw \rangle\bw-\langle \bw,\bw \rangle\bv)=\bw.
\end{eqnarray*}
\enD

We now denote 
\begin{eqnarray*}
A=\left|
\begin{array}{cc}
a_1 & b_1\\
a_2 & b_2
\end{array}
\right|,
B=\left|
\begin{array}{cc}
a_1 & c_1\\
a_2 & c_2
\end{array}
\right|,
C=\left|
\begin{array}{cc}
b_1 & c_1\\
b_2 & c_2
\end{array}
\right|.
\end{eqnarray*}
where $a_i,b_i,c_i:U \to \R,i=1,2$ are smooth functions defined as above. 
We remark that if  $(\bX,\bv,\bw)$ is a lightcone framed surface, then we have 
\begin{eqnarray*}
&a_1=b_1=0\ {\rm or}\ a_2=b_2=0, \\
&\bX_u\wedge\bX_v=-B\bv+C\bw.
\end{eqnarray*}
Without loss of generality, we assume that $ a_2=b_2=0$ in this paper.
This means that 
\begin{eqnarray*}
A=0,\ B=a_1c_2=-\alpha,\ C=b_1c_2=\beta 
\end{eqnarray*}
for all $(u,v) \in U$.

We note that a point $p\in U$ is said to be a lightlike point if the normal vector of $\bX$ at $p$ is lightlike; a point $p\in U$ is said to be a singular point if the normal vector of $\bX$ at $p$ vanishes.
As the normal vector of $\bX$ is $\bX_u\wedge\bX_v=-B\bv+C\bw=-a_1c_2\bv+b_1c_2\bw$. 
Then $p$ is a lightlike point of the lightcone framed surface $(\bX,\bv,\bw)$ if and only if $a_1(p)=0,b_1(p)\ne 0, c_2(p)\ne0$ or $a_1(p)\ne0,b_1(p)=0, c_2(p)\ne0$.
Without loss of generality, we may assume that $p$ is a lightlike point of the lightcone framed surface $(\bX,\bv,\bw)$ if and only if $a_1(p)=0$, $b_1(p)\ne0$ and $c_2(p)\ne0$.
The case for $b_1(p)=0$ is similar to the case for $a_1(p)=0$.
$p$ is a singular point of the lightcone framed surface $(\bX,\bv,\bw)$ if and only if $a_1(p)=b_1(p)=0$ or $c_2(p)=0$.

Let $p\in U$ be a singular point of the lightcone framed surface, if $a_1^2(p)+b_1^2(p)\ne0,c_2(p)=0$ (resp. $a_1^2(p)+b_1^2(p)=0,c_2(p)\ne0$), we call $p$ {\it the first singular point of the lightcone framed surface} (resp. {\it the second singular point of the lightcone framed surface}), and denote by
\begin{eqnarray*}
&&S_1(\bX)=\{a_1^2+b_1^2\ne0,c_2=0\},\begin{pmatrix}
 \bX_u\\\bX_v
\end{pmatrix}=\begin{pmatrix}
  a_1&  b_1& c_1\\
  0&  0& 0
\end{pmatrix}\begin{pmatrix}
 \bv\\
 \bw\\
 \bm
\end{pmatrix},\\
&&S_2(\bX)=\{a_1^2+b_1^2=0,c_2\ne0\},\begin{pmatrix}
 \bX_u\\\bX_v
\end{pmatrix}=\begin{pmatrix}
  0&  0& c_1\\
  0&  0& c_2
\end{pmatrix}\begin{pmatrix}
 \bv\\
 \bw\\
 \bm
\end{pmatrix}.
\end{eqnarray*}
In this paper, we investigate the lightcone framed surface with the first singular point.



\subsection{Curvatures of lightcone framed surfaces}

Here we calculate the Gaussian curvature and the Mean curvature of lightcone framed surfaces.

Since $\bX_v=c_2\bm=0$ at singular points of $\bX$, we denote the fundamental quantity by using $\bm$ instead of $\bX_v$. 
Then $\{\bX_u, \bm, \bn\}$ gives a frame at non-lightlike points of $\bX$, where $\bn=(\bX_u\wedge\bm) / ||\bX_u\wedge\bm||$.

We denote the first fundamental quantity of the lightcone framed surface $\bX$ with $\{\bX_u, \bm, \bn\}$ for each point $(u,v)\in U$ as followings:
\begin{eqnarray*}
&&\tilde{E}=\langle \bX_u,\bX_u \rangle=c^2_1-4a_1b_1,\\
&&\tilde{F}=\langle \bX_u,\bm \rangle=c_1,\\
&&\tilde{G}=\langle \bm,\bm \rangle=1.
\end{eqnarray*}

We define a function $\tilde{\lambda}:=\tilde{E}\tilde{G}-\tilde{F}^2=-4a_1b_1$ for each point $(u,v)\in U$.
We call $\tilde{\lambda}$ the {\it signed area density function} of the lightcone framed surface $(\bX,\bv,\bw)$.
Then a point p is a spacelike point of $\bX$ if $(a_1b_1)(p)<0$, lightlike point if $(a_1b_1)(p)=0$ and $(a_1^2+b_1^2)(p)\ne0$, timelike point if $(a_1b_1)(p)>0$ respectively. 
We denote the set of spacelike, lightlike and timelike points as $U_+$, $L(\bX)$ and $U_-$, respectively.

By a direct calculation, we know the normal vector of $\bX$ at non-lightlike points is
\begin{eqnarray*}
\bn=\frac{\bX_u\wedge\bm}{||\bX_u\wedge\bm||}=\frac{-a_1\bv+b_1\bw}{||-4a_1b_1||}=\frac{-a_1\bv+b_1\bw}{\sqrt{|\tilde{\lambda}|}}
\end{eqnarray*}
and
\begin{eqnarray*}
\bX_{uu}=(a_{1u}+a_1e_1+c_1f_1)\bv+(b_{1u}-b_1e_1+c_1g_1)\bw+(c_{1u}+2a_1g_1+2b_1f_1)\bm.
\end{eqnarray*}

Then we denote the second fundamental quantity of $\bX$ with $\{\bX_u, \bm, \bn\}$ at non-lightlike points as followings:
\begin{eqnarray*}
&&\tilde{L}=\langle \bX_{uu},\bn \rangle=\frac{2a_1(b_{1u}-b_1e_1+c_1g_1)-2b_1(a_{1u}+a_1e_1+c_1f_1)}{\sqrt{|\tilde{\lambda}|}},\\
&&\tilde{M}=\langle \bm_u,\bn \rangle=\frac{2(a_1g_1-b_1f_1)}{\sqrt{|\tilde{\lambda}|}},\\
&&\tilde{N}=\langle \bm_v,\bn \rangle=\frac{2(a_1g_2-b_1f_2)}{\sqrt{|\tilde{\lambda}|}}.
\end{eqnarray*}

In order to investigate the fundamental quantity of the lightcone framed surface with lightlike points, we define smooth functions

\begin{eqnarray*}
\hat{\bn}=\sqrt{|\tilde{\lambda}|}\bn, \hat{\tilde{L}}=\sqrt{|\tilde{\lambda}|}\tilde{L}, \hat{\tilde{M}}=\sqrt{|\tilde{\lambda}|}\tilde{M}, \hat{\tilde{N}}=\sqrt{|\tilde{\lambda}|}\tilde{N}.
\end{eqnarray*}

Then we have
\begin{eqnarray*}
&&\hat{\bn}=-a_1\bv+b_1\bw,\\
&&\hat{\tilde{L}}=2a_1(b_{1u}-b_1e_1+c_1g_1)-2b_1(a_{1u}+a_1e_1+c_1f_1),\\
&&\hat{\tilde{M}}=2(a_1g_1-b_1f_1),\\
&&\hat{\tilde{N}}=2(a_1g_2-b_1f_2).
\end{eqnarray*}
we denote $\hat{\bn}$ as the normal vector of the lightcone framed surface.

\begin{lemma}\label{relationship}
The relationships between the fundamental quantity $\tilde{E},\tilde{F},\tilde{G},\tilde{L},\tilde{M},\tilde{N}$ with $\{\bX_u,\bm,\bn\}$ and $E,F,G,L,M,N$ with $\{\bX_u,\bX_v,\bN\}$ are as follows, where $\bN=(\bX_u\wedge\bX_v) / || \bX_u\wedge\bX_v||$ :
\begin{eqnarray*}
E=\tilde{E}, F=c_2\tilde{F}, G=c_2^2\tilde{G}, \bN=\pm\bn, L=\tilde{L}, M=c_2\tilde{M}, N=c_2\tilde{N}.
\end{eqnarray*}
As the sign has no effect on the normal vector, we assume that $\bN=\bn$.
\end{lemma}
\demo
We know that
\begin{eqnarray*}
E=\langle \bX_u,\bX_u \rangle, F=\langle \bX_u,\bX_v \rangle, G=\langle \bX_v,\bX_v \rangle, L=\langle \bX_{uu},\bN \rangle, M=\langle \bX_{uv},\bN \rangle, N=\langle \bX_{vv},\bN \rangle,
\end{eqnarray*}
$$
\bN=\frac{\bX_u\wedge\bX_v}{||\bX_u\wedge\bX_v||}=\frac{c_2}{|c_2|}\frac{\bX_u\wedge\bm}{||\bX_u\wedge\bm||}=\pm\bn.
$$
We set $\bN=\bn$.
Thus by a direct calculation, we have
\begin{eqnarray*}
&&E=\langle \bX_u,\bX_u \rangle=\tilde{E}, F=c_2\langle \bX_u,\bm \rangle=c_2\tilde{F}, G=c_2^2\langle \bm,\bm \rangle=c_2^2\tilde{G}, L=\langle \bX_{uu},\bn \rangle=\tilde{L},\\
&&M=\langle \bX_{vu},\bn \rangle=\langle (c_2\bm)_u,\bn \rangle=\langle c_{2u}\bm+c_2\bm_u,\bn \rangle=c_2\langle \bm_u,\bn \rangle=c_2\tilde{M},\\
&&N=\langle \bX_{vv},\bn \rangle=\langle (c_2\bm)_v,\bn \rangle=\langle c_{2v}\bm+c_2\bm_v,\bn \rangle=c_2\langle \bm_v,\bn \rangle=c_2\tilde{N}.
\end{eqnarray*}
\enD

Now we consider the Gaussian curvature and the Mean curvature:
$$
K=\frac{LN-M^2}{|EG-F^2|}, H=\frac{LG-2MF+NE}{2|EG-F^2|}.
$$
By Lemma \ref{relationship}, we get the Gaussian curvature and the Mean curvature of $\bX$ at non-lightlike points:
\begin{eqnarray*}
&&K=\frac{c_2}{|c_2|}\frac{\tilde{L}\tilde{N}-c_2\tilde{M}^2}{|\tilde{E}\tilde{G}-\tilde{F}^2|}=\pm\frac{\tilde{L}\tilde{N}-c_2\tilde{M}^2}{|\tilde{E}\tilde{G}-\tilde{F}^2|},\\
&&H=\frac{c_2}{|c_2|}\frac{c_2\tilde{L}\tilde{G}-2c_2\tilde{M}\tilde{F}+\tilde{N}\tilde{E}}{2|\tilde{E}\tilde{G}-\tilde{F}^2|}=\pm\frac{c_2\tilde{L}\tilde{G}-2c_2\tilde{M}\tilde{F}+\tilde{N}\tilde{E}}{2|\tilde{E}\tilde{G}-\tilde{F}^2|}.
\end{eqnarray*}
where the sign of $K$ and $H$ are both $"+"$ if $c_2>0$. 
Then we have
\begin{eqnarray*}
K=&&\pm\frac{4a_1^2[g_2(b_{1u}-b_1e_1+c_1g_1)-c_2g_1^2]+4b_1^2[f_2(a_{1u}+a_1e_1+c_1f_1)-c_2f_1^2]}{|\tilde{\lambda}|^2}\\
&&\frac{-4a_1b_1[g_2(a_{1u}+a_1e_1+c_1f_1)+f_2(b_{1u}-b_1e_1+c_1g_1)-2c_2f_1g_1]}{|\tilde{\lambda}|^2},\\
H=&&\pm\frac{a_1[c_2(b_{1u}-b_1e_1-c_1g_1)+g_2(c_1^2-4a_1b_1)]-b_1[c_2(a_{1u}+a_1e_1-c_1f_1)+f_2(c_1^2-4a_1b_1)]}{|\tilde{\lambda}|^\frac{3}{2}}.
\end{eqnarray*}

In order to investigate the curvature of the lightcone framed surface with lightlike points, we define smooth functions 
$$
\hat{K}=|\tilde{\lambda}|^2K,\hat{H}=|\tilde{\lambda}|^\frac{3}{2}H,
$$
then we have
\begin{eqnarray*}
\hat{K}=&&\pm(4a_1^2[g_2(b_{1u}-b_1e_1+c_1g_1)-c_2g_1^2]+4b_1^2[f_2(a_{1u}+a_1e_1+c_1f_1)-c_2f_1^2]\\
&&-4a_1b_1[g_2(a_{1u}+a_1e_1+c_1f_1)+f_2(b_{1u}-b_1e_1+c_1g_1)-2c_2f_1g_1]),\\
\hat{H}=&&\pm (a_1[c_2(b_{1u}-b_1e_1-c_1g_1)+g_2(c_1^2-4a_1b_1)]- b_1[c_2(a_{1u}+a_1e_1-c_1f_1)+f_2(c_1^2-4a_1b_1)]).
\end{eqnarray*}

\begin{remark}\label{NLMNKH extend to U}
By the definition above, we can infer that 
\begin{eqnarray*}
&&\hat{\bn}=\bX_u \wedge \bm, \hat{\tilde{L}}=\langle \bX_{uu},\hat{\bn} \rangle, \hat{\tilde{M}}=\langle \bm_u,\hat{\bn} \rangle, 
\hat{\tilde{N}}=\langle \bm_v,\hat{\bn} \rangle,\\
&&\hat{K}=\pm(\hat{\tilde{L}}\hat{\tilde{N}}-c_2\hat{\tilde{M}}^2), 2\hat{H}=\pm(c_2\hat{\tilde{L}}\hat{\tilde{G}}-2c_2\hat{\tilde{M}}\hat{\tilde{F}}+\hat{\tilde{N}}\hat{\tilde{E}}).
\end{eqnarray*}
Then $\hat{\bn},\hat{\tilde{L}},\hat{\tilde{M}},\hat{\tilde{N}},\hat{K},\hat{H}$ can be extended to $U$.
It is obvious that $\tilde{\lambda}=-\langle \hat{\bn},\hat{\bn} \rangle$.
Morever $K=\hat{K}/|\tilde{\lambda}|^2$ and $H=\hat{H}/|\tilde{\lambda}|^\frac{3}{2}$ hold on non-lightlike points of the lightcone framed surface.
\end{remark}

\subsection{Behavior of lightlike points of lightcone framed surfaces}

In this part, we consider the properties of lightlike points of lightcone framed surfaces.
We first recall the notions of lightlike points of mixed type surfaces.
For more detail see \cite{Honda2}.

Let $f:U\to M^3$ be a mixed type surface, where $M^3$ is an oriented Lorentzian $3$-manifold.
For each point $p\in U$, we call the subspace
\begin{eqnarray*}
\mathcal{N}_p:=\{\nu \in T_pU; (ds^2)_p(\nu,\ba)=0,\forall \ba \in T_pU \}
\end{eqnarray*}
of $T_pU$ the {\it null space} at $p$.
Then $\mathcal{N}_p\ne\{\mathbf{0}\}$ if and only if $p$ is a lightlike point.
A nonzero element of $\mathcal{N}_p$ is called a {\it null vector} at $p$.
We remark that $dim\mathcal{N}_p=2$ does not occur for any lightlike point $p$, as $df_p(T_pU)$ is a $2$-dimensional degenerate subspace of the Lorentzian product space $T_{f(p)}M^3$.
Then we have
\begin{eqnarray*}
dim\mathcal{N}_p=1
\end{eqnarray*}
for each lightlike point $p$ of a mixed type surface $f$.
Then a smooth vector $\eta$ defined on a neighborhood $U$ of a lightlike point $p$ is called a {\it null vector field} if $\eta_q\in T_qU$ is a null vector at each $q\in LD\cap U$, where $LD$ is the set of lightlike points of $f$.
This means that $df(\eta)$ is lightlike at a lightlike point $p$.

On a local coordinate neighborhood $(U;u,v)$, set $f_u:=df(\partial_u),f_v:=df(\partial_v)$, where $\partial_u:=\partial/\partial_u,\partial_v:=\partial/\partial_v$.
Then 
\begin{eqnarray*}
ds^2=Edu^2+2Fdudv+Gdv^2,
\end{eqnarray*}
where $E:=\langle f_u,f_u \rangle,F:=\langle f_u,f_v \rangle,G:=\langle f_v,f_v \rangle$.
We set a function $\lambda:=EG-F^2$.
A lightlike point $p\in LD$ is called {\it non-degenerate} if $d\lambda(p)\ne0$, {\it degenerate} if $d\lambda(p)=0$.
By the implicit function theorem, there exists a regular curve $\gamma(t)(|t|<\varepsilon)$ in $U$ such that $p=\gamma(0)$ and $Image(\gamma)=LD$ hold on a neighborhood of $p$.
We call $\gamma(t)$ a {\it characteristic curve}.
For a null vector field $\eta$ defined on a neighborhood of $p$, the restriction $\eta(t):=\eta_{\gamma(t)}\in T_{\gamma(t)}U$ is called a {\it null vector field} along $\gamma(t)$.

Let $p\in U$ be a non-degenerate lightlike point, if $\gamma'(0)$ and $\eta(0)$ are linearly independent (resp. linearly dependent), we call $p$ the {\it lightlike point of the first kind} (resp. the {\it lightlike point of the second kind}).

\begin{lemma}\label{lightlike points of X}
Let $p\in U$ be a lightlike point of the lightcone framed surface $(\bX,\bv,\bw)$. 
We said that $p$ is a non-degenerate (resp. degenerate) lightlike point of $\bX$ if $a_1(p)=0, (a_{1u}^2+a_{1v}^2)(p)\ne0$ (resp. $a_1(p)=a_{1u}(p)=a_{1v}(p)=0$).
\end{lemma}
\demo
As the discriminant function $\tilde{\lambda}(p)=(-4a_1b_1)(p)$, then
\begin{eqnarray*}
&&\tilde{\lambda}_u(p)=(-4(a_{1u}b_1+a_1b_{1u}))(p),\\
&&\tilde{\lambda}_v(p)=(-4(a_{1v}b_1+a_1b_{1v}))(p).
\end{eqnarray*}
As $a_1(p)=0$ with $b_1(p)\ne0, c_2(p)\ne0$, we have 
\begin{eqnarray*}
&&\tilde{\lambda}_u(p)=(-4a_{1u}b_1)(p),\\
&&\tilde{\lambda}_v(p)=(-4a_{1v}b_1)(p).
\end{eqnarray*}
Thus $d\tilde{\lambda}(p)\ne0$ if and only if $\tilde{\lambda}_u^2(p)+\tilde{\lambda}_v^2(p)\ne0$ if and only if $a_{1u}^2(p)+a_{1v}^2(p)\ne0$;
$d\tilde{\lambda}(p)=0$ if and only if $\tilde{\lambda}_u(p)=\tilde{\lambda}_v(p)=0$ if and only if $a_{1u}(p)=a_{1v}(p)=0$.
\enD

\begin{lemma}
Let $p\in U$ be a lightlike point of the lightcone framed surface $(\bX,\bv,\bw)$.
Then the null vector at $p$ is 
\begin{eqnarray*}
\eta=c_2\partial_u-c_1\partial_v.
\end{eqnarray*}
\end{lemma}
\demo
Let the null vector at $p$ be $\eta=\eta_1\partial_u+\eta_2\partial_v$.
Then 
\begin{eqnarray*}
(d\bX(\eta))(p)=(\eta_1\bX_u+\eta_2\bX_v)(p)=(\eta_1(a_1\bv+b_1\bw+c_1\bm)+\eta_2c_2\bm)(p)=(\eta_1a_1\bv+\eta_1b_1\bw+(\eta_1c_1+\eta_2c_2)\bm)(p).
\end{eqnarray*}
As the vector $d\bX(\eta)$ at $p$ is lightlike, we have
\begin{eqnarray}\label{null vector of X}
(\langle d\bX(\eta),d\bX(\eta) \rangle)(p)=(\eta_1^2(c_1^2-4a_1b_1)+\eta_2^2c_2^2+2\eta_1\eta_2c_1c_2)(p)=0.
\end{eqnarray}
It is obvious that $\eta_1=c_2$ and $\eta_2=-c_1$ satisfy \eqref{null vector of X}.
So we get the null vector $\eta=c_2\partial_u-c_1\partial_v$ at a lightlike point $p$.
\enD

Next we investigate the properties of lightlike points of lightcone framed surfaces.
It is known that the induced metric of a mixed type surface is a signature-changing metric and the lightlike points may be regarded as singular points of such metrics.
So we consider lightlike points of mixed type surfaces in the way similar to the case of singular points of wave fronts(\cite{Honda2}).
The criteria theorem of singularities of fronts is as follows:

\begin{theorem}[\cite{Teramoto1}]\label{criterion-singularities}
Let $f:(U,p)\to\R^3$ be a front germ, let $\nu$ be a unit normal to $f$ and let $p$ be a corank one singular point; thus rank $df_p=1$.\\
$(1)$ Suppose that $p$ is a non-degenerate singular point of $f$, then
\par 
$(i)$ $f$ at $p$ is $\mathcal{A}$-equivalent to a cuspidal edge if and only if $\eta\lambda(p)\ne0$;\\
\par 
$(ii)$ $f$ at $p$ is $\mathcal{A}$-equivalent to a swallowtail if and only if $\eta\lambda(p)=0$ and $\eta\eta\lambda(p)\ne0$;\\
\par 
$(iii)$ $f$ at $p$ is $\mathcal{A}$-equivalent to a cuspidal butterfly if and only if $\eta\lambda(p)=\eta\eta\lambda(p)=0$ and $\eta\eta\eta\lambda(p)\ne0$;\\
$(2)$ Suppose that $p$ is a degenerate singular point of $f$, then
\par 
$(i)$ $f$ at $p$ is $\mathcal{A}$-equivalent to a cuspidal lips if and only if $det\mathcal{H}_{\lambda}(p)>0$;\\
\par 
$(ii)$ $f$ at $p$ is $\mathcal{A}$-equivalent to a cuspidal beaks if and only if $\eta\eta\lambda(p)\ne0$ and $det\mathcal{H}_{\lambda}(p)<0$;\\
Here $\lambda$ is the signed area density function of $f$, $\eta$ is the null vector field and $\mathcal{H}_{\lambda}$ is the Hessian matrix of $\lambda$.
\end{theorem}

Now we give the classification of lightlike points of lightcone framed surfaces.

\begin{theorem}
Let $(\bX,\bv,\bw):U \to \R_1^3 \times \Delta_4$ be a lightcone framed surface and $p$ be a lightlike point of $\bX$. \\
$(1)$ If $\bX$ at $p$ is non-degeneracy, we have that $a_1(p)=0,(a_{1u}^2+a_{1v}^2)(p)\ne0$. Then\\
\par
$(a)$ $\bX$ at $p$ is $\mathcal{A}$-equivalent to the cuspidal-edge-kind if and only if $(a_{1v}c_1-a_{1u}c_2)(p)\ne0$; \\
\par
$(b)$ $\bX$ at $p$ is $\mathcal{A}$-equivalent to the swallowtail-kind if and only if \\
\begin{eqnarray*}
(a_{1v}c_1-a_{1u}c_2)(p)=0
\end{eqnarray*}
and
\begin{eqnarray*}
(c_2(a_{1v}c_{1u}-a_{1uu}c_2-a_{1u}c_{2u}+a_{1uv}c_1)-c_1(a_{1vv}c_1+a_{1v}c_{1v}-a_{1u}c_{2v}-a_{1uv}c_2))(p)\ne0;
\end{eqnarray*}
\par
$(c)$ $\bX$ at $p$ is $\mathcal{A}$-equivalent to the cuspidal-butterfly-kind if and only if \\
\begin{eqnarray*}
(a_{1v}c_1-a_{1u}c_2)(p)=(c_2(a_{1v}c_{1u}-a_{1uu}c_2-a_{1u}c_{2u}+a_{1uv}c_1)-c_1(a_{1vv}c_1+a_{1v}c_{1v}-a_{1u}c_{2v}-a_{1uv}c_2))(p)=0
\end{eqnarray*}
and
\begin{eqnarray*}
&&((c_2c_{2u}-c_1c_{2v})(a_{1v}c_{1u}-a_{1uu}c_2-a_{1u}c_{2u}+a_{1uv}c_1)\\
&&+(c_1c_{1v}-c_2c_{1u})(a_{1vv}c_1+a_{1v}c_{1v}-a_{1u}c_{2v}-a_{1uv}c_2)\\
&&+c_2^2(2a_{1uv}c_{1u}+a_{1v}c_{1uu}-a_{1uuu}c_2-2a_{1uu}c_{2u}-a_{1u}c_{2uu}+a_{1uuv}c_1)\\
&&+c_1^2(a_{1vvv}c_1+2a_{1vv}c_{1v}+a_{1v}c_{1vv}-2a_{1uv}c_{2v}-a_{1u}c_{2vv}-a_{1uvv}c_2)\\
&&-2c_1c_2(a_{1uvv}c_1+a_{1vv}c_{1u}+a_{1uv}c_{1v}+a_{1v}c_{1uv}-a_{1uu}c_{2v}-a_{1u}c_{2uv}-a_{1uuv}c_2-a_{1uv}c_{2u}))(p)\ne0.
\end{eqnarray*}\\
$(2)$If $\bX$ at $p$ is degeneracy, we have that $a_1(p)=a_{1u}(p)=a_{1v}(p)=0$. Then\\
\par
$(a)$ $\bX$ at $p$ is $\mathcal{A}$-equivalent to the cuspidal-lips-kind if and only if $(a_{1uu}a_{1vv}-a_{1uv}^2)(p)>0$; \\
\par
$(b)$ $\bX$ at $p$ is $\mathcal{A}$-equivalent to the cuspidal-beaks-kind if and only if 
$$
(a_{1uu}a_{1vv}-a_{1uv}^2)(p)<0
$$
and
$$
(2a_{1uv}c_1c_2-a_{1uu}c_2^2-a_{1vv}c_1^2)(p)\ne0.
$$
\end{theorem}
\demo
We know that the null vector $\eta=c_2\partial_u-c_1\partial_v$ at lightlike points of $\bX$ and the signed area density function of $\bX$ is $\tilde{\lambda}=-4a_1b_1$.
Thus by a direct calculation, we have
$$
\tilde{\lambda}_u=-4(a_{1u}b_1+a_1b_{1u}),\tilde{\lambda}_v=-4(a_{1v}b_1+a_1b_{1v}).
$$

First we prove the assertion $(1)$. 
If $\bX$ at $p$ is non-degeneracy, we have $a_1(p)=0,(a_{1u}^2+a_{1v}^2)(p)\ne0$.
So 
\begin{eqnarray*}
\eta\tilde{\lambda}=c_2\tilde{\lambda}_u-c_1\tilde{\lambda}_v &&=-4c_2(a_{1u}b_1+a_1b_{1u})+4c_1(a_{1v}b_1+a_1b_{1v})\\
&&=4a_1(b_{1v}c_1-b_{1u}c_2)+4b_1(a_{1v}c_1-a_{1u}c_2)
\end{eqnarray*}

As $a_1(p)=0,b_1(p)\ne0$, we have 
$$
\eta\tilde{\lambda}(p)=(4b_1(a_{1v}c_1-a_{1u}c_2))(p),
$$
Thus $\eta\tilde{\lambda}(p)\ne0$ if and only if 
$$
(a_{1v}c_1-a_{1u}c_2)(p)\ne0.
$$
Therefore by Theorem \ref{criterion-singularities}, we obtain the conclusion $(a)$.

For the canclusion $(b)$, $\eta\tilde{\lambda}(p)=0$ if and only if 
$$
(a_{1v}c_1-a_{1u}c_2)(p)=0.
$$
Differentiate $\eta\tilde{\lambda}$ with respect to $u$ and $v$ respectively, we have
\begin{eqnarray*}
(\eta\tilde{\lambda})_u=&4a_{1u}(b_{1v}c_1-b_{1u}c_2)+4a_1(b_{1uv}c_1+b_{1v}c_{1u}-b_{1uu}c_1-b_{1u}c_{2u})\\
&+4b_{1u}(a_{1v}c_1-a_{1u}c_2)+4b_1(a_{1uv}c_1+a_{1v}c_{1u}-a_{1uu}c_2-a_{1u}c_{2u},
\end{eqnarray*}
\begin{eqnarray*}
(\eta\tilde{\lambda})_v=&4a_{1v}(b_{1v}c_1-b_{1u}c_2)+4a_1(b_{1vv}c_1+b_{1v}c_{1v}-b_{1uv}c_2-b_{1u}c_{2v})\\
&+4b_{1v}(a_{1v}c_1-a_{1u}c_2)+4b_1(a_{1vv}c_1+a_{1v}c_{1v}-a_{1uv}c_2-a_{1u}c_{2v}).
\end{eqnarray*}
Then 
\begin{eqnarray*}
\eta\eta\tilde{\lambda}=&&c_2(\eta\tilde{\lambda})_u-c_1(\eta\tilde{\lambda})_v\\
=&&4a_{1u}c_2(b_{1v}c_1-b_{1u}c_2)+4a_1c_2(b_{1uv}c_1+b_{1v}c_{1u}-b_{1uu}c_1-b_{1u}c_{2u})\\
&&+4b_{1u}c_2(a_{1v}c_1-a_{1u}c_2)+4b_1c_2(a_{1uv}c_1+a_{1v}c_{1u}-a_{1uu}c_2-a_{1u}c_{2u}\\
&&-4a_{1v}c_1(b_{1v}c_1-b_{1u}c_2)-4a_1c_1(b_{1vv}c_1+b_{1v}c_{1v}-b_{1uv}c_2-b_{1u}c_{2v})\\
&&-4b_{1v}c_1(a_{1v}c_1-a_{1u}c_2)-4b_1c_1(a_{1vv}c_1+a_{1v}c_{1v}-a_{1uv}c_2-a_{1u}c_{2v})\\
=&&8(a_{1u}c_2-a_{1v}c_1)(b_{1v}c_1-b_{1u}c_2)\\
&&+4a_1[c_2(b_{1v}c_{1u}-b_{1uu}c_2-b_{1u}c_{2u}+b_{1uv}c_1)-c_1(b_{1vv}c_1+b_{1v}c_{1v}-b_{1u}c_{2v}-b_{1uv}c_2)]\\
&&+4b_1[c_2(a_{1v}c_{1u}-a_{1uu}c_2-a_{1u}c_{2u}+a_{1uv}c_1)-c_1(a_{1vv}c_1+a_{1v}c_{1v}-a_{1u}c_{2v}-a_{1uv}c_2)].
\end{eqnarray*}
As $a_1(p)=0$, $b_1(p)\ne0$ and $(a_{1v}c_1-a_{1u}c_2)(p)=0$, we have
\begin{eqnarray*}
\eta\eta\tilde{\lambda}(p)=(4b_1[c_2(a_{1v}c_{1u}-a_{1uu}c_2-a_{1u}c_{2u}+a_{1uv}c_1)-c_1(a_{1vv}c_1+a_{1v}c_{1v}-a_{1u}c_{2v}-a_{1uv}c_2)])(p),
\end{eqnarray*}
Thus $\eta\eta\tilde{\lambda}(p)\ne0$ if and only if 
\begin{eqnarray*}
(c_2(a_{1v}c_{1u}-a_{1uu}c_2-a_{1u}c_{2u}+a_{1uv}c_1)-c_1(a_{1vv}c_1+a_{1v}c_{1v}-a_{1u}c_{2v}-a_{1uv}c_2))(p)\ne0
\end{eqnarray*}
Then by Theorem \ref{criterion-singularities}, we obtain the conclusion $(b)$.

For the canclusion $(c)$, $\eta\tilde{\lambda}(p)=0$ if and only if 
\begin{eqnarray}\label{3-9-1}
(a_{1v}c_1-a_{1u}c_2)(p)=0,
\end{eqnarray}
and $\eta\eta\tilde{\lambda}(p)=0$ if and only if 
\begin{eqnarray}\label{3-9-2}
(c_2(a_{1v}c_{1u}-a_{1uu}c_2-a_{1u}c_{2u}+a_{1uv}c_1)-c_1(a_{1vv}c_1+a_{1v}c_{1v}-a_{1u}c_{2v}-a_{1uv}c_2))(p)=0.
\end{eqnarray}
By a direct calculation and the conditions: $a_1(p)=0$, \eqref{3-9-1} and \eqref{3-9-2}, we have
\begin{eqnarray*}
\eta\eta\eta\tilde{\lambda}(p)=&&(4b_1[(c_2c_{2u}-c_1c_{2v})(a_{1v}c_{1u}-a_{1uu}c_2-a_{1u}c_{2u}+a_{1uv}c_1)\\
&&+(c_1c_{1v}-c_2c_{1u})(a_{1vv}c_1+a_{1v}c_{1v}-a_{1u}c_{2v}-a_{1uv}c_2)\\
&&+c_2^2(2a_{1uv}c_{1u}+a_{1v}c_{1uu}-a_{1uuu}c_2-2a_{1uu}c_{2u}-a_{1u}c_{2uu}+a_{1uuv}c_1)\\
&&+c_1^2(a_{1vvv}c_1+2a_{1vv}c_{1v}+a_{1v}c_{1vv}-2a_{1uv}c_{2v}-a_{1u}c_{2vv}-a_{1uvv}c_2)\\
&&-2c_1c_2(a_{1uvv}c_1+a_{1vv}c_{1u}+a_{1uv}c_{1v}+a_{1v}c_{1uv}-a_{1uu}c_{2v}-a_{1u}c_{2uv}-a_{1uuv}c_2-a_{1uv}c_{2u})])(p).
\end{eqnarray*}
Since $b_1(p)\ne0$, then $\eta\eta\eta\tilde{\lambda}(p)\ne0$ if and only if 
\begin{eqnarray*}
&&((c_2c_{2u}-c_1c_{2v})(a_{1v}c_{1u}-a_{1uu}c_2-a_{1u}c_{2u}+a_{1uv}c_1)\\
&&+(c_1c_{1v}-c_2c_{1u})(a_{1vv}c_1+a_{1v}c_{1v}-a_{1u}c_{2v}-a_{1uv}c_2)\\
&&+c_2^2(2a_{1uv}c_{1u}+a_{1v}c_{1uu}-a_{1uuu}c_2-2a_{1uu}c_{2u}-a_{1u}c_{2uu}+a_{1uuv}c_1)\\
&&+c_1^2(a_{1vvv}c_1+2a_{1vv}c_{1v}+a_{1v}c_{1vv}-2a_{1uv}c_{2v}-a_{1u}c_{2vv}-a_{1uvv}c_2)\\
&&-2c_1c_2(a_{1uvv}c_1+a_{1vv}c_{1u}+a_{1uv}c_{1v}+a_{1v}c_{1uv}-a_{1uu}c_{2v}-a_{1u}c_{2uv}-a_{1uuv}c_2-a_{1uv}c_{2u}))(p)\ne0.
\end{eqnarray*}
Thus by Theorem \ref{criterion-singularities}, we obtain the canclusion $(c)$.

Now we prove the assertion$(2)$.
If $\bX$ at $p$ is degeneracy, we have $a_1(p)=a_{1u}(p)=a_{1v}(p)=0$.
By a direct calculation, we have
\begin{eqnarray*}
&&\tilde{\lambda}_{uu}(p)=(-4(a_{1uu}b_1+2a_{1u}b_{1u}+a_1b_{1uu}))(p)=(-4a_{1uu}b_1)(p),\\
&&\tilde{\lambda}_{uv}(p)=(-4(a_{1uv}b_1+a_{1u}b_{1v}+a_{1v}b_{1u}+a_1b_{1uv}))(p)=(-4a_{1uv}b_1)(p),\\
&&\tilde{\lambda}_{vv}(p)=(-4(a_{1vv}b_1+2a_{1v}b_{1v}+a_1b_{1vv}))(p)=(-4a_{1vv}b_1)(p).
\end{eqnarray*}
The Hessian matrix of $\tilde{\lambda}$ is 
\begin{eqnarray*}
det\mathcal{H}_{\tilde{\lambda}}(p)=\begin{vmatrix}
  \tilde{\lambda}_{uu}&\tilde{\lambda}_{uv} \\
  \tilde{\lambda}_{uv}&\tilde{\lambda}_{vv}
\end{vmatrix}(p)=\begin{vmatrix}
  -4a_{1uu}b_1&-4a_{1uv}b_1 \\
  -4a_{1uv}b_1&-4a_{1vv}b_1
\end{vmatrix}(p)=(16b_1^2(a_{1uu}a_{1vv}-a_{1uv}^2))(p).
\end{eqnarray*}
So $det\mathcal{H}_{\tilde{\lambda}}(p)>0$ if and only if 
$$
(a_{1uu}a_{1vv}-a_{1uv}^2)(p)>0,
$$
$det\mathcal{H}_{\tilde{\lambda}}(p)<0$ if and only if 
$$
(a_{1uu}a_{1vv}-a_{1uv}^2)(p)<0.
$$
By the condition $a_1(p)=a_{1u}(p)=a_{1v}(p)=0$, we get
$$
\eta\eta\tilde{\lambda}(p)=(4b_1(2a_{1uv}c_1c_2-a_{1uu}c_2^2-a_{1vv}c_1^2)(p).
$$
Then $\eta\eta\tilde{\lambda}(p)\ne0$ if and only if 
$$
(2a_{1uv}c_1c_2-a_{1uu}c_2^2-a_{1vv}c_1^2)(p)\ne0.
$$
Thus by Theorem \ref{criterion-singularities}, we have the result. 
\enD

\section{Focal surfaces of lightcone framed surfaces}

In this section, we investigate the differential geometry of focal surfaces of lightcone framed surfaces.

\subsection{Principal curvature and principal direction of lightcone framed surfaces}

We consider the principal curvature and principal direction of the lightcone framed surface $\bX$ by using Weingarten transformation.
We set the Weingarten map $W$ at non-lightlike points of $\bX$ as follows:
$$
W\begin{pmatrix}
 \bX_u\\\bX_v
\end{pmatrix}=-\begin{pmatrix}
 \bN_u\\\bN_v
\end{pmatrix},
$$
Then multiply both sides by $\begin{pmatrix}
  \bX_u&\bX_v
\end{pmatrix}$,we get
\begin{eqnarray*}
W=&&\begin{pmatrix}
  L&M \\
  M&N
\end{pmatrix}
\begin{pmatrix}
  E&F \\
  F&G
\end{pmatrix}^{-1}
=\begin{pmatrix}
  \tilde{L}&c_2\tilde{M} \\
  c_2\tilde{M}&c_2\tilde{N}
\end{pmatrix}
\begin{pmatrix}
  \tilde{E}&c_2\tilde{F} \\
  c_2\tilde{F}&c_2^2\tilde{G}
\end{pmatrix}^{-1}\\
=&&\begin{pmatrix}
  \tilde{L}&c_2\tilde{M} \\
  c_2\tilde{M}&c_2\tilde{N}
\end{pmatrix}
\begin{pmatrix}
  \frac{\tilde{G}}{\tilde{E}\tilde{G}-\tilde{F}^2}&\frac{-\tilde{F}}{c_2(\tilde{E}\tilde{G}-\tilde{F}^2)} \\
  \frac{-\tilde{F}}{c_2(\tilde{E}\tilde{G}-\tilde{F}^2)}&\frac{\tilde{E}}{c_2^2(\tilde{E}\tilde{G}-\tilde{F}^2)}
\end{pmatrix}\\
=&&\begin{pmatrix}
  \frac{\tilde{L}\tilde{G}-\tilde{M}\tilde{F}}{\tilde{E}\tilde{G}-\tilde{F}^2}&\frac{\tilde{M}\tilde{E}-\tilde{L}\tilde{F}}{c_2(\tilde{E}\tilde{G}-\tilde{F}^2)} \\
  \frac{c_2\tilde{M}\tilde{G}-\tilde{N}\tilde{F}}{\tilde{E}\tilde{G}-\tilde{F}^2}&\frac{\tilde{N}\tilde{E}-c_2\tilde{M}\tilde{F}}{c_2(\tilde{E}\tilde{G}-\tilde{F}^2)}
\end{pmatrix}
\end{eqnarray*}
We know the fact that the eigenvalue of the matrix $W$ is the principal curvature of the surface $\bX$. 
We denote $\kappa_i (i=1,2)$ as the principal curvature of surface $\bX$ at non-lightlike part, then

\begin{eqnarray*}
|\kappa_iI-W|=&&\begin{vmatrix}
  \kappa_i-\frac{\tilde{L}\tilde{G}-\tilde{M}\tilde{F}}{\tilde{E}\tilde{G}-\tilde{F}^2}&\frac{\tilde{L}\tilde{F}-\tilde{M}\tilde{E}}{c_2(\tilde{E}\tilde{G}-\tilde{F}^2)} \\
  \frac{\tilde{N}\tilde{F}-c_2\tilde{M}\tilde{G}}{\tilde{E}\tilde{G}-\tilde{F}^2}&\kappa_i-\frac{\tilde{N}\tilde{E}-c_2\tilde{M}\tilde{F}}{c_2(\tilde{E}\tilde{G}-\tilde{F}^2)}
\end{vmatrix}\\
=&&\kappa_i^2-\frac{c_2\hat{\tilde{L}}\hat{\tilde{G}}-2c_2\hat{\tilde{M}}\hat{\tilde{F}}+\hat{\tilde{N}}\hat{\tilde{E}}}{c_2(\tilde{E}\tilde{G}-\tilde{F}^2)}\kappa_i+\frac{(\tilde{L}\tilde{N}-c_2\tilde{M}^2)(\tilde{E}\tilde{G}-\tilde{F}^2)}{c_2(\tilde{E}\tilde{G}-\tilde{F}^2)^2}=0,(i=1,2).
\end{eqnarray*}
that is
$$
c_2\tilde{\lambda}\kappa_i^2-2|\tilde{\lambda}|H\kappa_i+|\tilde{\lambda}|K=0.
$$
So the principal curvature at non-lightlike and non-singular points is
$$
\kappa_i=\frac{H\pm \sqrt{H^2-c_2K}}{c_2},(i=1,2), \ when \ \tilde{\lambda}>0,
$$
or
$$
\kappa_i=\frac{-H\pm \sqrt{H^2+c_2K}}{c_2},(i=1,2), \ when \ \tilde{\lambda}<0.
$$

As $\hat{K}=|\tilde{\lambda}|^2K,\hat{H}=|\tilde{\lambda}|^\frac{3}{2}H$, we have 
\begin{eqnarray}\label{ki1}
c_2\tilde{\lambda}|\tilde{\lambda}|\kappa_i^2-2|\tilde{\lambda}|^\frac{1}{2}\hat{H}\kappa_i+\hat{K}=0.
\end{eqnarray}
then
$$
\kappa_i=\frac{\hat{H}\pm\sqrt{\hat{H}^2-c_2\tilde{\lambda}\hat{K}}}{c_2\tilde{\lambda}|\tilde{\lambda}|^\frac{1}{2}},
$$

In order to consider the principal curvature of lightcone framed surfaces with lightlike points, we set $\hat{\kappa_i}=\tilde{\lambda}|\tilde{\lambda}|^\frac{1}{2}\kappa_i,(i=1,2)$, 
then 
$$
\hat{\kappa_i}=\frac{\hat{H}\pm\sqrt{\hat{H}^2-c_2\tilde{\lambda}\hat{K}}}{c_2},(i=1,2).
$$

Note that $\hat{\kappa_i}$ is the principal curvature at non-singular points.
When $p$ is the singular point of $\bX$, then by (\ref{ki1}) we have 
$$
\kappa_i(p)=\frac{\hat{K}(p)}{2|\tilde{\lambda}(p)|^\frac{1}{2}\hat{H}(p)}.
$$
We set $\hat{\kappa_i}=|\tilde{\lambda}|^\frac{1}{2}\kappa_i$, 
then the principal curvature at singular points is
$$
\hat{\kappa_i}(p)=\frac{\hat{K}(p)}{2\hat{H}(p)}.
$$

Let $\bV_i=(\bV_{i1},\bV_{i2})(i=1,2)$ be the principal vector of the lightcone framed surface at non-lightlike part, then as $(II-\kappa_i I)\bV_i=\mathbf{0}, (i=1,2)$, we have

$$
\left [\begin{pmatrix}
  L&M \\
  M&N
\end{pmatrix}-\kappa_i\begin{pmatrix}
  E&F \\
  F&G
\end{pmatrix}\right ]\begin{pmatrix}
 \bV_{i1}\\\bV_{i2}
\end{pmatrix}=\begin{pmatrix}
 0\\0
\end{pmatrix},(i=1,2),
$$
that is
$$
\begin{pmatrix}
  L-\kappa_i E&M-\kappa_i F \\
  M-\kappa_i F&N-\kappa_i G
\end{pmatrix}\begin{pmatrix}
 \bV_{i1}\\\bV_{i2}
\end{pmatrix}=\begin{pmatrix}
 0\\0
\end{pmatrix},(i=1,2).
$$
Then
$$
\bV_i=(N-\kappa_i G,-M+\kappa_i F)=(\tilde{N}-c_2\kappa_i\tilde{G},-\tilde{M}+\kappa_i\tilde{F}),(i=1,2).
$$

As $\hat{\kappa_i}=\tilde{\lambda}|\tilde{\lambda}|^\frac{1}{2}\kappa_i(i=1,2),\hat{\tilde{M}}=\sqrt{|\tilde{\lambda}|}\tilde{M,}\hat{\tilde{N}}=\sqrt{|\tilde{\lambda}|}\tilde{N}$, we have 
$$
\bV_i=\frac{1}{|\lambda|^\frac{1}{2}}(\hat{\tilde{N}}-\frac{c_2}{\tilde{\lambda}}\hat{\kappa_i}\tilde{G},-\hat{\tilde{M}}+\frac{1}{\tilde{\lambda}}\hat{\kappa_i}\tilde{F}),(i=1,2).
$$
Set $\hat{\bV_i}=|\lambda|^\frac{1}{2}\bV_i,(i=1,2)$, then 
$$
\hat{\bV_i}=(\tilde{\lambda}\hat{\tilde{N}}-c_2\hat{\kappa_i}\tilde{G},-\tilde{\lambda}\hat{\tilde{M}}+\hat{\kappa_i}\tilde{F}),(i=1,2).
$$

\begin{remark}
$\hat{\kappa_i}$ and $\hat{\bV_i}(i=1,2)$ can be extended to $U$.
\end{remark}

\subsection{Focal surfaces of lightcone framed surfaces}

For points on $L(\bX)$, the normal line to $\bX$ is lightlike and is tangent to $\bX$.
The focal set of $\bX$ is not defined at points on $L(\bX)$.
So we consider the caustic of $\bX$ as the bifurcation set of the family of distance-squared function on $\bX$(cf. \cite{Tari1}).
Note that the caustic coincides with the focal set of $\bX\setminus L(\bX)$ and can extend the focal set to the $L(\bX)$.
Thus we define focal surfaces of lightcone framed surfaces by caustics as follows.

The family of distance-squared functions on $\bX$ is given by 
$$
d^2:\bX \times \R_1^3 \to \R,(p,\nu)\mapsto d^2(p,\nu)=\langle p-\nu,p-\nu \rangle
$$
We denote by 
$$
\Sigma(d^2)=\{\xi=((u,v),\nu)\in U\times\R_1^3 \ | \ d_u^2(\xi)=d_v^2(\xi)=0\},
$$
and by
\begin{eqnarray*}
Bif(d^2)=\{\nu\in\R_1^3 \ | \ \exists((u,v),\nu)\in\Sigma(d^2), rank(\mathcal{H}(d^2)(u,v))<2\}.
\end{eqnarray*}

The set $Bif(d^2)$ is the local stratum of the {\it bifurcation set} of the family $d^2$, that is, it is the set of points $\nu\in\R_1^3$ for which there exists $(u,v)\in U$ such that $d^2$ has a degenerate local singularity at $(u,v)$.
Then we can get a parametrization of the caustic as follows. 
Since $d^2(u,v)=\langle \bX-\nu,\bX-\nu\rangle$, $d^2$ is singular at $p\in U$ if and only if $\langle \bX-\nu,\bX_u\rangle=\langle \bX-\nu,\bX_v\rangle=0$ at $p$, if and only if $\bX-\nu=\mu_0\bX_u\times \bX_v=\mu_0c_2\bX_u\wedge\bm=\mu\hat{\bn}$ for some scalar $\mu_0$ and $\mu$.
Note that the condition includes the case when $p$ is on the $L(\bX)$ where $\bX_u\times \bX_v$ is parallel to the unique lightlike direction at $p$.

The singularity of $d^2$ at $p$ is degenerate if and only if $\bX-\nu=\mu\hat{\bn}$ for some scalar $\mu$ and $(d^2_{uv})^2(p)-d^2_{uu}(p)d^2_{vv}(p)=0$,where
\begin{eqnarray*}
&&\frac{1}{2}d^2_{uu}=\langle \bX_u,\bX_u \rangle+\langle \bX-\nu,\bX_{uu} \rangle
=\langle \bX_u,\bX_u \rangle+\mu\langle \hat{\bn},\bX_{uu}\rangle\\
&&\frac{1}{2}d^2_{uv}=\langle \bX_u,\bX_v \rangle+\langle \bX-\nu,\bX_{uv} \rangle
=c_2\langle \bX_u,\bm \rangle+\mu\langle \hat{\bn},\bX_{uv}\rangle,\\
&&\frac{1}{2}d^2_{vv}=\langle \bX_v,\bX_v \rangle+\langle \bX-\nu,\bX_{vv} \rangle
=c_2^2\langle \bm,\bm \rangle+\mu\langle \hat{\bn},\bX_{vv}\rangle.
\end{eqnarray*}

Since 
\begin{eqnarray*}
&&\langle \hat{\bn},\bX_{uv}\rangle=\langle \hat{\bn}, (c_2\bm)_u \rangle=\langle \hat{\bn}, c_{2u}\bm+c_2\bm_u \rangle
=c_2\langle \hat{\bn}, \bm_u \rangle,\\
&&\langle \hat{\bn},\bX_{vv}\rangle=\langle \hat{\bn}, (c_2\bm)_v \rangle=\langle \hat{\bn}, c_{2v}\bm+c_2\bm_v \rangle
=c_2\langle \hat{\bn}, \bm_v \rangle,
\end{eqnarray*}
then by Remark \ref{NLMNKH extend to U} we have 
\begin{eqnarray*}
&&\frac{1}{2}d^2_{uu}=\tilde{E}+\mu\hat{\tilde{L}},\\
&&\frac{1}{2}d^2_{uv}=c_2\tilde{F}+\mu c_2\hat{\tilde{M}},\\
&&\frac{1}{2}d^2_{vv}=c_2^2\tilde{G}+\mu c_2\hat{\tilde{N}}.
\end{eqnarray*}
Therefore the singularity of $d^2$ at $p$ is degenerate if and only if $\bX-\nu=\mu\hat{\bn}$ and
$$
(\hat{\tilde{L}}\hat{\tilde{N}}-c_2\hat{\tilde{M}}^2)\mu^2+(c_2\hat{\tilde{L}}\tilde{G}-2c_2\hat{\tilde{M}}\tilde{F}+\hat{\tilde{N}}\tilde{E})\mu+c_2(\tilde{E}\tilde{G}-\tilde{F}^2)=0,
$$
that is 
$$
\hat{K}\mu^2+2\hat{H}\mu+c_2\tilde{\lambda}=0,
$$
So
$$
\mu=-\frac{\hat{H}\pm\sqrt{\hat{H}^2-c_2\tilde{\lambda}\hat{K}}}{\hat{K}}.
$$

Then we obtain focal surfaces of lightcone framed surfaces as follows:
\begin{definition}
Let $(\bX,\bv,\bw):U \to \R_1^3\times\Delta_4$ be a lightcone framed surface.
We call $\bF:U \to \R_1^3, \bF=\bX-\mu\hat{\bn}$ the {\it focal surface} of the lightcone framed surface, where 
$$
\mu=-\frac{\hat{H}\pm\sqrt{\hat{H}^2-c_2\tilde{\lambda}\hat{K}}}{\hat{K}}.
$$
\end{definition}

\begin{remark}\label{mu at lightlike and singular points}
If $(u,v)$ is a lightlike or singular point of $\bX$, we have $\mu=0$ or $\mu=-\frac{2\hat{H}}{\hat{K}}$.
\end{remark}

We now consider the differential geometry of focal surfaces of lightcone framed surfaces in $\R_1^3$.
By a direct calculation, we have:
\begin{eqnarray*}
&&\bF_u=[a_1+\mu(a_{1u}+a_1e_1)+\mu_ua_1]\bv+[b_1-\mu(b_{1u}-b_1e_1)-\mu_ub_1]\bw+[c_1+2\mu(a_1g_1-b_1f_1)]\bm,\\
&&\bF_v=[\mu(a_{1v}+a_1e_2)+\mu_va_1]\bv-[\mu(b_{1v}-b_1e_2)+\mu_vb_1]\bw+[c_2+2\mu(a_1g_2-b_1f_2)]\bm,\\
\end{eqnarray*}

Then we can give the following formula:
$$
\begin{pmatrix}
 \bF_u\\\bF_v
\end{pmatrix}:=\begin{pmatrix}
  \bar{a}_1&  \bar{b}_1& \bar{c}_1\\
  \bar{a}_2&  \bar{b}_2& \bar{c}_2
\end{pmatrix}\begin{pmatrix}
 \bv\\
 \bw\\
 \bm
\end{pmatrix},
$$
where
\begin{eqnarray}\label{invariants of F}
&&\bar{a}_1=a_1(1+\mu e_1+\mu_u)+\mu a_{1u},\\
&&\bar{b}_1=b_1(1+\mu e_1-\mu_u)-\mu b_{1u},\\
&&\bar{c}_1=c_1+2\mu (a_1g_1-b_1f_1),\\
&&\bar{a}_2=a_1(\mu e_2+\mu_v)+\mu a_{1v},\\
&&\bar{b}_2=-b_1(-\mu e_2+\mu_v)-\mu b_{1v},\\
&&\bar{c}_2=c_2+2\mu(a_1g_2-b_1f_2).
\end{eqnarray}

\begin{proposition}
Under the above settings, the focal surface $\bF$ is a lightcone framed surface with the lightcone frame $\{\bv,\bw,\bm\}$ if and only if $a_1(\mu e_2+\mu_v)+\mu a_{1v}=a_1(\mu e_2+\mu_v)+\mu a_{1v}=0$ for each $(u,v)\in U$.
In this case, the formula above can be written as:
$$
\begin{pmatrix}
 \bF_u\\\bF_v
\end{pmatrix}:=\begin{pmatrix}
  \bar{a}_1& \bar{b}_1 & \bar{c}_1\\
  0&  0& \bar{c}_2
\end{pmatrix}\begin{pmatrix}
 \bv\\
 \bw\\
 \bm
\end{pmatrix}.
$$
\end{proposition}

According to the fact that if the normal vector of a surface is lightlike at a point $p\in U$, $p$ is a lightlike point of the surface; if the normal vector of a surface vanishes at a point $p\in U$, $p$ is a singular point of the surface, we have the following results:
\begin{proposition}
Let $(\bX,\bv,\bw):U \to \R_1^3\times\Delta_4$ be a lightcone framed surface and $\bF:U\to\R_1^3$ be the focal surface of $\bX$. 
Then $(u,v)$ is a singular point of $\bF$ if and only if $\bar{a}_1=\bar{b}_1=0$ or $\bar{c}_2=0$.
And we denote by
\begin{eqnarray*}
&&S_1(\bF)=\{\bar{a}_1^2+\bar{b}_1^2\ne0,\bar{c}_2=0\},\begin{pmatrix}
 \bF_u\\\bF_v
\end{pmatrix}=\begin{pmatrix}
  \bar{a}_1&  \bar{b}_1& \bar{c}_1\\
  0&  0& 0
\end{pmatrix}\begin{pmatrix}
 \bv\\
 \bw\\
 \bm
\end{pmatrix},\\
&&S_2(\bF)=\{\bar{a}_1^2+\bar{b}_1^2=0,\bar{c}_2\ne0\},\begin{pmatrix}
 \bF_u\\\bF_v
\end{pmatrix}=\begin{pmatrix}
  0&  0& \bar{c}_1\\
  0&  0& \bar{c}_2
\end{pmatrix}\begin{pmatrix}
 \bv\\
 \bw\\
 \bm
\end{pmatrix}.
\end{eqnarray*}
\end{proposition}

\begin{proposition}\label{lightlike points of F}
Let $(\bX,\bv,\bw):U \to \R_1^3\times\Delta_4$ be a lightcone framed surface and $\bF:U\to\R_1^3$ be the focal surface of $\bX$.
Then $p$ is a lightlike point of $\bF$ if and only if and $\bar{a}_1(p)=0,\bar{b}_1(p)\ne0,\bar{c}_2(p)\ne0$ or $\bar{a}_1(p)\ne0,\bar{b}_1(p)=0,\bar{c}_2(p)\ne0$.
Without loss of generality, we suppose that $p$ is a lightlike point of $\bF$ if and only if $\bar{a}_1(p)=0,\bar{b}_1(p)\ne0$ and $\bar{c}_2(p)\ne0$.
The case for $\bar{b}_1(p)=0$ is similar to the case for $\bar{a}_1(p)=0$.
\end{proposition}

The relationship of lightlike and singular points between the focal surface and the original surface is as follows:
\begin{theorem}
Let $(\bX,\bv,\bw):U \to \R_1^3\times\Delta_4$ be a lightcone framed surface and $\bF:U\to\R_1^3$ be the focal surface of $\bX$.
If $p$ is a denegeracy lightlike point of $\bX$, then $p$ is a lightlike point of $\bF$.
\end{theorem}
\demo
It is known that $p$ is a degeneracy lightlike point of $\bX$ if and only if $b_1(p)\ne0, c_2(p)\ne0, a_1(p)=a_{1u}(p)=a_{1v}(p)=0$.
Then by \eqref{invariants of F}, we get $\bar{a}_1(p)=0$.
So by Proposition \ref{lightlike points of F}, $p$ is a lightlike point of $F$.
\enD

\begin{theorem}
Let $(\bX,\bv,\bw):U \to \R_1^3\times\Delta_4$ be a lightcone framed surface and $\bF:U\to\R_1^3$ be the focal surface of $\bX$.
If $(a_1g_2)(p)=(b_1f_2)(p)$, then $p$ is a singular point of $\bX$ if and only if $p$ is a singular point of $\bF$.
\end{theorem}
\demo
When $p$ is a singular point of $\bX$, we have $c_2(p)=0$.
As $\bar{c}_2(p)=(c_2+2\mu(a_1g_2-b_1f_2))(p)$ and $a_1g_2=b_1f_2$, we have $\bar{c}_2(p)=0$.
Then we get the calculation.
\enD

\subsection{Curvature of focal surfaces of lightcone framed surfaces}

In this part, we consider the Gaussian curvature and Mean curvature of focal surfaces of lightcone framed surfaces in $\R_1^3$.

Since $\bF_v=\bar{c}_2\bm=0$ at singular points of $\bF$, we denote the fundamental quantity by using $\bm$ instead of $\bF_v$. 
Then $\{\bF_u, \bm, \bar{\bn}\}$ gives a frame at non-lightlike points of $\bF$, where $\bar{\bn}=(\bF_u\wedge\bm) / ||\bF_u\wedge\bm||$.

We denote the first fundamental quantity of the focal surface $\bF$ of the lightcone framed surface $\bX$ with $\{\bF_u, \bm, \bar{\bn}\}$ for each point $(u,v)\in U$ as followings:

The first fundamental quantity of $\bF$ for each point $(u,v)\in U$ is as follows:
\begin{eqnarray*}
&&\tilde{\bar{E}}=\langle \bF_u,\bF_u \rangle=\bar{c}^2_1-4\bar{a}_1\bar{b}_1,\\
&&\tilde{\bar{F}}=\langle \bF_u,\bm \rangle=\bar{c}_1,\\
&&\tilde{\bar{G}}=\langle \bm,\bm \rangle=1.
\end{eqnarray*}

We define a function $\tilde{\bar{\lambda}}:=\tilde{\bar{E}}\tilde{\bar{G}}-\tilde{\bar{F}}^2=-4\bar{a}_1\bar{b}_1$ for each point $(u,v)\in U$.
We call $\tilde{\bar{\lambda}}$ the {\it signed area density function} of the focal surface of the lightcone framed surface.
Then a point p is a spacelike point of $\bF$ if $(\bar{a}_1\bar{b}_1)(p)<0$, lightlike point if $(\bar{a}_1\bar{b}_1)(p)=0$ and $(\bar{a}_1^2+\bar{b}_1^2)(p)\ne0$, timelike point if $(\bar{a}_1\bar{b}_1)(p)>0$ respectively. 

By a direct calculation, we know the normal vector of $\bF$ at non-lightlike points is
\begin{eqnarray*}
\bar{\bn}=\frac{\bF_u\wedge\bm}{||\bF_u\wedge\bm||}=\frac{-\bar{a}_1\bv+\bar{b}_1\bw}{||-4\bar{a}_1\bar{b}_1||}=\frac{-\bar{a}_1\bv+\bar{b}_1\bw}{\sqrt{|\tilde{\bar{\lambda}}|}}
\end{eqnarray*}

Now we denote the second fundamental quantity of $\bF$ with $\{\bF_u, \bm, \bar{\bn}\}$ at non-lightlike points as followings:
\begin{eqnarray*}
&&\tilde{\bar{L}}=\langle \bF_{uu},\bar{\bn}\rangle 
=\frac{2\bar{a}_1(\bar{b}_{1u}-\bar{b}_1e_1+\bar{c}_1g_1)-2\bar{b}_1(\bar{a}_{1u}+\bar{a}_1e_1+\bar{c}_1f_1)}{\sqrt{|\tilde{\bar{\lambda}}|}},\\
&&\tilde{\bar{M}}=\langle \bm_u,\bar{\bn}\rangle
=\frac{2(\bar{a}_1g_1-\bar{b}_1f_1)}{\sqrt{|\tilde{\bar{\lambda}}|}},\\
&&\tilde{\bar{N}}=\langle \bm_v,\bar{\bn}\rangle
=\frac{2(\bar{a}_1g_2-\bar{b}_1f_2)}{\sqrt{|\tilde{\bar{\lambda}}|}}.
\end{eqnarray*}

In order to investigate the fundamental quantity of the focal surface of the lightcone framed surface with lightlike points, we define smooth functions
\begin{eqnarray*}
\hat{\bar{\bn}}=\sqrt{|\tilde{\bar{\lambda}}|}\bar{\bn}, \hat{\tilde{\bar{L}}}=\sqrt{|\tilde{\bar{\lambda}}|}\tilde{\bar{L}}, \hat{\tilde{\bar{M}}}=\sqrt{|\tilde{\bar{\lambda}}|}\tilde{\bar{M}}, \hat{\tilde{\bar{N}}}=\sqrt{|\tilde{\bar{\lambda}}|}\tilde{\bar{N}}.
\end{eqnarray*}

Then we have
\begin{eqnarray*}
&&\hat{\bar{\bn}}=-\bar{a}_1\bv+\bar{b}_1\bw,\\
&&\hat{\tilde{\bar{L}}}=2\bar{a}_1(\bar{b}_{1u}-\bar{b}_1e_1+\bar{c}_1g_1)-2\bar{b}_1(\bar{a}_{1u}+\bar{a}_1e_1+\bar{c}_1f_1),\\
&&\hat{\tilde{\bar{M}}}=2(\bar{a}_1g_1-\bar{b}_1f_1),\\
&&\hat{\tilde{\bar{N}}}=2(\bar{a}_1g_2-\bar{b}_1f_2).
\end{eqnarray*}

\begin{lemma}\label{relationship1}
The relationships between the fundamental quantity $\tilde{\bar{E}},\tilde{\bar{F}},\tilde{\bar{G}},\tilde{\bar{L}},\tilde{\bar{M}},\tilde{\bar{N}}$ with $\{\bF_u,\bm,\bar{\bn}\}$ and $\bar{E},\bar{F},\bar{G},\bar{L},\bar{M},\bar{N}$ with $\{\bF_u,\bF_v,\bar{\bN}\}$ are as follows, where $\bar{\bN}=(\bF_u\wedge\bF_v) / || \bF_u\wedge\bF_v||$ :
\begin{eqnarray*}
\bar{E}=\tilde{\bar{E}}, \bar{F}=\bar{c}_2\tilde{\bar{F}}, \bar{G}=\bar{c}_2^2\tilde{\bar{G}}, \bar{\bN}=\pm\bar{\bn}, \bar{L}=\tilde{\bar{L}}, \bar{M}=\bar{c}_2\tilde{\bar{M}}, \bar{N}=\bar{c}_2\tilde{\bar{N}}.
\end{eqnarray*}
As the sign has no effect on the normal vector, we assume that $\bar{\bN}=\bar{\bn}$.
\end{lemma}

Now we consider the Gaussian curvature and the Mean curvature of the focal surface:
By Lemma \ref{relationship1}, we get the Gaussian curvature and the Mean curvature of $\bF$ at non-lightlike points:
\begin{eqnarray*}
&&\bar{K}=\pm\frac{\tilde{\bar{L}}\tilde{\bar{N}}-\bar{c}_2\tilde{\bar{M}}^2}{|\tilde{\bar{E}}\tilde{\bar{G}}-\tilde{\bar{F}}^2|},\\
&&\bar{H}=\pm\frac{\bar{c}_2\tilde{\bar{L}}\tilde{\bar{G}}-2\bar{c}_2\tilde{\bar{M}}\tilde{\bar{F}}+\tilde{\bar{N}}\tilde{\bar{E}}}{2|\tilde{\bar{E}}\tilde{\bar{G}}-\tilde{\bar{F}}^2|}.
\end{eqnarray*}
where the sign of $\bar{K}$ and $\bar{H}$ are both $"+"$ if $\bar{c}_2>0$.
Then we have
\begin{eqnarray*}
\bar{K}=&&\pm\frac{4\bar{a}_1^2[g_2(\bar{b}_{1u}-\bar{b}_1e_1+\bar{c}_1g_1)-\bar{c}_2g_1^2]+4\bar{b}_1^2[f_2(\bar{a}_{1u}+\bar{a}_1e_1+\bar{c}_1f_1)-\bar{c}_2f_1^2]}{|\tilde{\bar{\lambda}}|^2}\\
&&\frac{-4\bar{a}_1\bar{b}_1[g_2(\bar{a}_{1u}+\bar{a}_1e_1+\bar{c}_1f_1)+f_2(\bar{b}_{1u}-\bar{b}_1e_1+\bar{c}_1g_1)-2\bar{c}_2f_1g_1]}{|\tilde{\bar{\lambda}}|^2},\\
\bar{H}=&&\pm\frac{\bar{a}_1[\bar{c}_2(\bar{b}_{1u}-\bar{b}_1e_1-\bar{c}_1g_1)+g_2(\bar{c}_1^2-4\bar{a}_1\bar{b}_1)]-\bar{b}_1[\bar{c}_2(\bar{a}_{1u}+\bar{a}_1e_1-\bar{c}_1f_1)+f_2(\bar{c}_1^2-4\bar{a}_1\bar{b}_1)]}{|\tilde{\bar{\lambda}}|^\frac{3}{2}}.
\end{eqnarray*}

In order to investigate the curvature of the focal surface of lightcone framed surfaces with lightlike points, we define smooth functions 
$$
\hat{\bar{K}}=|\tilde{\bar{\lambda}}|^2\bar{K},\hat{\bar{H}}=|\tilde{\bar{\lambda}}|^\frac{3}{2}\bar{H},
$$
Then we have
\begin{eqnarray*}
\hat{\bar{K}}=&&\pm(4\bar{a}_1^2[g_2(\bar{b}_{1u}-\bar{b}_1e_1+\bar{c}_1g_1)-\bar{c}_2g_1^2]+4\bar{b}_1^2[f_2(\bar{a}_{1u}+\bar{a}_1e_1+\bar{c}_1f_1)-\bar{c}_2f_1^2]\\
&&-4\bar{a}_1\bar{b}_1[g_2(\bar{a}_{1u}+\bar{a}_1e_1+\bar{c}_1f_1)+f_2(\bar{b}_{1u}-\bar{b}_1e_1+\bar{c}_1g_1)-2\bar{c}_2f_1g_1]),\\
\hat{\bar{H}}=&&\pm(\bar{a}_1[\bar{c}_2(\bar{b}_{1u}-\bar{b}_1e_1-\bar{c}_1g_1)+g_2(\bar{c}_1^2-4\bar{a}_1\bar{b}_1)]-\bar{b}_1[\bar{c}_2(\bar{a}_{1u}+\bar{a}_1e_1-\bar{c}_1f_1)-f_2(\bar{c}_1^2-4\bar{a}_1\bar{b}_1)]).
\end{eqnarray*}

\begin{remark}
By the definition above, we can infer that
\begin{eqnarray*}
&&\hat{\bar{\bn}}=\bF_u\wedge\bm,
\hat{\tilde{\bar{L}}}=\langle \bF_{uu},\hat{\bar{\bn}} \rangle,
\hat{\tilde{\bar{M}}}=\langle \bm_u,\hat{\bar{\bn}} \rangle,
\hat{\tilde{\bar{N}}}=\langle \bm_v,\hat{\bar{\bn}} \rangle,\\
&&\hat{\bar{K}}=\pm(\hat{\tilde{\bar{L}}}\hat{\tilde{\bar{N}}}-\bar{c}_2\hat{\tilde{\bar{M}}}^2), 2\hat{\bar{H}}=\pm(\bar{c}_2\hat{\tilde{\bar{L}}}\hat{\tilde{\bar{G}}}-2\bar{c}_2\hat{\tilde{\bar{M}}}\hat{\tilde{\bar{F}}}+\hat{\tilde{\bar{N}}}\hat{\tilde{\bar{E}}}).
\end{eqnarray*}
Then $\hat{\bar{\bn}},\hat{\tilde{\bar{L}}},\hat{\tilde{\bar{M}}},\hat{\tilde{\bar{N}}},\hat{\bar{K}},\hat{\bar{H}}$ can be extended to $U$.
\end{remark}

\subsection{Differential geometry property of focal surfaces of lightcone framed surfaces}

\begin{remark}
As $K=\hat{K}/|\tilde{\lambda}|^2,H=\hat{H}/|\tilde{\lambda}|^\frac{3}{2}$, if $p$ is a lightlike point of $\bX$ and $\hat{K}(p)\ne0$, then $K$ is unbounded at $p$, and if $p$ is a lightlike point of $\bX$ and $\hat{H}(p)\ne0$, then $H$ is unbounded at $p$. 
Note that $K$ and $H$ cannot be both bounded at $p$ (for details see \cite{Honda2}).            
\end{remark}

We will characterize the properties of the Mean curvature and the Gaussian curvature of the lightcone framed surface $\bX$ and the focal surface $\bF$ of $\bX$ as follows.

\begin{theorem}
Let $(\bX,\bv,\bw):U \to \R_1^3\times\Delta_4$ be a lightcone framed surface and $\bF:U\to\R_1^3$ be a focal surface of $\bX$.
Assume that the Mean curvature $H$ of $\bX$ is bounded on $U$ and $p$ is a lightlike point of $\bX$.
Then $\hat{H}(p)=0$ and $p$ is a lightlike point of $\bF$.
Morever $\hat{\bar{H}}(p)=0$ if and only if $((\mu_u-1)\mu_ua_{1u})(p)=0$.
Specially, if the lightlike point $p$ of $\bX$ is degeneracy, then $\hat{\bar{H}}(p)=0$.
\end{theorem}
\demo
Since $H=\hat{H}/|\tilde{\lambda}|^\frac{3}{2}$ and $\tilde{\lambda}(p)=0$, we can infer that $\hat{H}(p)=0$.
In fact, if $\hat{H}(p)\ne0$ and $\tilde{\lambda}(p)=0$, we get $H$ is unbounded.
This contradicts the assumption.
By Remark \ref{lightlike points of X}, $p$ is a lightlike point of $\bX$ if and only if $a_1(p)=0$.
Then 
$$
\hat{H}(p)=(\mp b_1[c_2(a_{1u}-c_1f_1)+c_1^2f_2])(p)=0,
$$
Since $b_1(p)\ne0$, we obtain that 
$$
(c_2(a_{1u}-c_1f_1)+c_1^2f_2)(p)=0.
$$

By Remark \ref{mu at lightlike and singular points}, $\mu(p)=0 \ or \ -2\hat{H}(p)/\hat{K}(p)$.
As the fact that the Gaussian curvature and Mean curvature cannot be both bounded at lightlike points of mixed type surfaces(cf. \cite{Honda2}), we have that $K(p)$ is undounded.
Then since $K=\hat{K}/|\tilde{\lambda}|^2$ and $\tilde{\lambda}(p)=0$, we have $\hat{K}(p)=0$ or $\hat{K}(p)=K_0\ne0$ where $K_0$ is a constant in $\R$.

When $\hat{K}(p)=0$, as 
$$
\hat{H}(p)=(\mp b_1[c_2(a_{1u}-c_1f_1)+c_1^2f_2])(p),\hat{K}(p)=(\pm 4b_1^2[f_2(a_{1u}+c_1f_1)-c_2f_1^2])(p),
$$
we see that the degree of the denominator $\hat{K}(p)$ of $\mu(p)=(-2\hat{H}/\hat{K})(p)$ is greater than that of the numerator $\hat{H}(p)$.
This means that the denominator $\hat{K}(p)$ of $\mu(p)$ approaches zero faster than the numerator $\hat{H}(p)$.
Thus $\mu(p)=0$ for this case.
So we get $\mu(p)=0$.
Then 
\begin{eqnarray*}
&&\bar{a}_1(p)=(a_1(1+\mu e_1+\mu_u)+\mu a_{1u})(p)=0,\\
&&\bar{b}_1(p)=(b_1(1+\mu e_1-\mu_u)-\mu b_{1u})(p)=(b_1(1-\mu_u ))(p),\\
&&\bar{c}_1(p)=(c_1+2\mu (a_1g_1-b_1f_1))(p)=c_1(p),\\
&&\bar{c}_2(p)=(c_2+2\mu(a_1g_2-b_1f_2))(p)=c_2(p)\ne0,\\
&&\bar{a}_{1u}(p)=(a_{1u}(1+\mu e_1+\mu_u)+\mu_u a_{1u})(p)=(a_{1u}(1+2\mu_u))(p),
\end{eqnarray*}
By Proposition \ref{lightlike points of F}, $p$ is a lightlike point of $\bF$.

When $\hat{K}(p)=K_0\ne0$, we get $\mu(p)=0$ or $\mu(p)=-2\hat{H}(p)/\hat{K}(p)=0/K_0=0$.

By a direct calculation, 
\begin{eqnarray*}
\hat{\bar{H}}(p)=&&(\mp\bar{b}_1[\bar{c}_2(\bar{a}_{1u}-\bar{c}_1f_1)+\bar{c}_1^2f_2))(p)\\
=&&(\pm(\mu_u-1)b_1[c_2(a_{1u}+2\mu_u a_{1u}-c_1f_1)+c_1^2f_2]\\
=&&(\pm2(\mu_u-1)\mu_ua_{1u}b_1c_2)(p),
\end{eqnarray*}
So $\hat{\bar{H}}(p)=0$ if and only if 
$$
((\mu_u-1)\mu_ua_{1u})(p)=0.
$$

Specially, if the lightlike point $p$ of $\bX$ is degeneracy, we must have $\hat{\bar{H}}(p)=0$.
\enD

\begin{remark}
If $((\mu_u-1)\mu_ua_{1u})(p)\ne0$, $\bar{H}(p)$ is unbounded. 
\end{remark}

\begin{theorem}
Let $(\bX,\bv,\bw):U \to \R_1^3\times\Delta_4$ be a lightcone framed surface and $\bF:U\to\R_1^3$ be the focal surface of $\bX$.
Assume that the Gaussian curvature $K$ of $\bX$ is bounded on $U$ and $p$ is a lightlike point of $\bX$.
Then $\hat{K}(p)=0$ and one of the followings is true: \\
$(1)$ 
$\mu(p)=0$, $p$ is a lightlike point of $\bF$ and $\hat{\bar{K}}(p)=0$ if and only if $((1-\mu_u)\mu_ua_{1u}f_2)(p)=0$.
Specially, if the lightlike point $p$ of $\bX$ is degeneracy, $\hat{\bar{K}}(p)=0$.\\
$(2)$
$\mu(p)$ is unbounded and $p$ is a non-lightlike point of $\bF$ if and only if $a_{1u}(p)\ne0, \bar{b}_1(p)\ne0$ and $f_2(p)\ne0$. 
\end{theorem}
\demo
Since $K=\hat{K}/|\tilde{\lambda}|^2$ and $\tilde{\lambda}(p)=0$, we can infer that $\hat{K}(p)=0$.
In fact, if $\hat{K}(p)\ne0$ and $\tilde{\lambda}(p)=0$, we get $K$ is unbounded.
This contradicts the assumption.
By Remark \ref{lightlike points of X}, $p$ is a lightlike point of $\bX$ if and only if $a_1(p)=0$.
Then 
$$
\hat{K}(p)=\pm(4b_1^2[f_2(a_{1u}+c_1f_1)-c_2f_1^2])(p)=0,
$$
Since $(b_1)(p)\ne0$, we obtain that 
$$
(f_2(a_{1u}+c_1f_1)-c_2f_1^2)(p)=0,
$$

By Remark \ref{mu at lightlike and singular points}, $\mu(p)=0 \ or \ -2\hat{H}(p)/\hat{K}(p)$.
As the fact that the Gaussian curvature and Mean curvature can not be both bounded at lightlike points of mixed type surfaces, $H(p)$ is undounded.
Then since $H=\hat{H}/|\tilde{\lambda}|^\frac{3}{2}$ and $\tilde{\lambda}(p)=0$, we have $\hat{H}(p)=0$ or $\hat{H}(p)=H_0\ne0$ where $H_0$ is a constant in $\R$.

When $\hat{H}(p)=0$, as 
$$
\hat{H}(p)=(\mp b_1[c_2(a_{1u}-c_1f_1)+c_1^2f_2])(p),\hat{K}(p)=(\pm4b_1^2[f_2(a_{1u}+c_1f_1)-c_2f_1^2])(p),
$$
we see that the degree of the denominator $\hat{K}(p)$ of $\mu(p)=(-2\hat{H}/\hat{K})(p)$ is greater than that of the numerator $\hat{H}(p)$.
This means that the denominator $\hat{K}(p)$ of $\mu(p)$ approaches zero faster than the numerator $\hat{H}(p)$.
Thus $\mu(p)=0$ for this case.
Then 
\begin{eqnarray*}
&&\bar{a}_1(p)=(a_1(1+\mu e_1+\mu_u)+\mu a_{1u})(p)=0,\\
&&\bar{b}_1(p)=(b_1(1+\mu e_1-\mu_u)-\mu b_{1u})(p)=(b_1(1-\mu_u))(p),\\
&&\bar{c}_1(p)=(c_1+2\mu (a_1g_1-b_1f_1))(p)=c_1(p),\\
&&\bar{c}_2(p)=(c_2+2\mu(a_1g_2-b_1f_2))(p)=c_2(p)\ne0,\\
&&\bar{a}_{1u}(p)=(a_{1u}(1+\mu_u)+\mu_ua_{1u})(p)=(a_{1u}(1+2\mu_u))(p),
\end{eqnarray*}
By Proposition \ref{lightlike points of F}, $p$ is a lightlike point of $\bF$.
By a direct calculation, 
\begin{eqnarray*}
\hat{\bar{K}}(p)=&&(\pm4\bar{b}_1^2[f_2(\bar{a}_{1u}+\bar{c}_1f_1)-\bar{c}_2f_1^2])(p)\\
=&&(\pm4(1-\mu_u)^2b_1^2[(a_{1u}f_2(2\mu_u+1)+c_1f_1f_2-c_2f_1^2])(p)\\
=&&(\pm8(1-\mu_u)^2\mu_ua_{1u}b_1^2f_2)(p),
\end{eqnarray*}
So $\hat{\bar{K}}(p)=0$ if and only if 
$$
((1-\mu_u)\mu_ua_{1u}f_2)(p)=0.
$$
Specially, if the lightlike point $p$ of $\bX$ is degeneracy, we must have $\hat{\bar{K}}(p)=0$.

When $\hat{H}(p)=H_0\ne0$, we get $\mu(p)$ is unbounded.
Then
\begin{eqnarray*}
&&\bar{a}_1(p)=(a_1(1+\mu e_1+\mu_u)+\mu a_{1u})(p)=(\mu a_{1u}c_2)(p),\\
&&\bar{b}_1(p)=(b_1(1+\mu e_1-\mu_u)-\mu b_{1u})(p),\\
&&\bar{c}_2(p)=(c_2+2\mu(a_1g_2-b_1f_2))(p)=(c_2-2\mu b_1f_2)(p),
\end{eqnarray*}
we know that $\bar{a}_1(p)$ is unbounded if $a_{1u}(p)\ne0$, and $\bar{c}_2(p)$ is unbounded if $f_2(p)\ne0$.
So $\bar{a}_1(p)\ne0$ and $\bar{c}_2(p)\ne0$ if $f_2(p)\ne0$.
Therefore by Proposition \ref{lightlike points of F}, $p$ is a non-lightlike point of $\bF$ if $a_{1u}(p)\ne0, \bar{b}_1(p)\ne0$ and $f_2(p)\ne0$.
\enD

\begin{theorem}(cf. \cite{Honda1})
Let $f:U \to (M_3,g)$ a real analytic immersion into an oriented real analytic Lorentzian manifold $(M_3,g)$.
Denote by $U_+$(resp. $U_-$) the set of points where $f$ is spacelike(resp. timelike).
Suppose that $U_+,U_-$ are both non-empty and the Mean curvature $H$ on $U_+\cup U_-$ is bounded.
Then for each $p \in \overline{U_+}\cap \overline{U_-}$, there exists a sequence $\{p_n\}_{n=1,2,3,...}$ in $U_+$(resp. $U_-$) converging to $p$ such that $\lim_{n \to \infty} H(p_n)=0$, where $\overline{U_+},\overline{U_-}$ are the closures of $U_+,U_-$ in $U$.
\end{theorem}

\begin{corollary}\label{the Mean curvature of F} 
Let $(\bX,\bv,\bw):U \to \R_1^3\times\Delta_4$ be a lightcone framed surface and $\bF:U\to\R_1^3$ be the focal surface of $\bX$.
Suppose that the Mean curvature $H$ of $\bX$ and $\bar{H}$ of $\bF$ are both bounded at each lightlike point $p$ of $\bX$.
Then there exists a sequence $\{p_n\}_{n=1,2,3,...}$ in $U_+$(resp. $U_-$) converging to $p$ such that 
$$
\lim_{n \to \infty} \bar{H}(p_n)=0.
$$
\end{corollary}
\demo
Let $q_1,q_2\in V\subset U$ satisfy $\tilde{\bar{\lambda}}(q_1)>0,\tilde{\bar{\lambda}}(q_2)<0$ and $\gamma(s),s\in[0,2\pi]$ be a smooth curve in $U$ satisfied $\gamma(0)=q_1,\gamma(2\pi)=q_2$.
Then by Fourier series expansion of $\gamma(s)$, there exists a sufficiently large positive integer $N$, the real analytic curve defined by $\gamma_N(s)$ satisfies
\begin{eqnarray}\label{lambda-gamma-1}
\tilde{\bar{\lambda}}(\gamma_N(0))>0,\tilde{\bar{\lambda}}(\gamma_N(2\pi))<0.
\end{eqnarray}
We set 
$$
\check{\bar{\lambda}}(s):=\tilde{\bar{\lambda}}(\gamma_N(s)),s\in[0,2\pi]
$$
Since $\check{\bar{\lambda}}(s)$ is a real analytic function defined on $[0,2\pi]$, the set of zeros of the function $\check{\bar{\lambda}}(s)$ consists of a finite set of points
$$
0<s_1<s_2<\cdot\cdot\cdot<s_n<2\pi.
$$
By \eqref{lambda-gamma-1}, we can choose the number $j$ such that the sign of $\check{\bar{\lambda}}(s)$ changes from positive to negative at $s=s_j$.
So there exists a positive integer $m$ such that
$$
\lim_{s \to s_j}\frac{\check{\bar{\lambda}}(s)}{(s-s_j)^m}=\bar{\lambda}_0\ne0,
$$
where $\bar{\lambda}_0$ is a non-zero real number.
As $\check{\bar{\lambda}}(s)$ changes sign at $s=s_j$, the number $m$ is odd.
By Remark \ref{NLMNKH extend to U}, $\hat{H}$ can be extended to $U$.
We set
$$
\check{\hat{\bar{H}}}(s):=\hat{\bar{H}}(\gamma_N(s)),s\in[0,2\pi],
$$
then 
$$
\bar{H}(\gamma_N(s)):=\frac{\check{\hat{\bar{H}}}(s)}{|\check{\bar{\lambda}}(s)|^\frac{3}{2}},s\ne s_1,s_2,\cdot\cdot\cdot,s_n.
$$
As $\bar{H}(p)$ is bounded at each $p\in L(\bX)$, we get $\hat{\bar{H}}(s_j)=0$.
Since $\hat{\bar{H}}(s)=0$ is a real analytic function, there exists a positive integer $l$ such that 
$$
\lim_{s \to s_j} \frac{\check{\hat{\bar{H}}}(s)}{(s-s_j)^l}=\bar{H}_0\ne0,
$$
where $\bar{H}_0$ is a non-zero real number.
Then it holds that
\begin{eqnarray*}
\bar{H}(\gamma_N(s))=\frac{\check{\hat{\bar{H}}}(s)}{(s-s_j)^l} \frac{(s-s_j)^{\frac{3}{2}m}}{|\check{\bar{\lambda}}(s)|^\frac{3}{2}}(s-s_j)^{l-\frac{3}{2}m},
\end{eqnarray*}
then
\begin{eqnarray*}
\lim_{s \to s_j} |s-s_j|^{\frac{3}{2}m-l}|\bar{H}(\gamma_N(s))|=\frac{|\bar{H}_0|}{|\bar{\lambda}_0|^\frac{3}{2}}\ne0.
\end{eqnarray*}
As $\bar{H}(p)$ is bounded, we have $2l-3m\ge0$.
Since $m$ is odd, we have $2l-3m>0$.
Then 
$$
\lim_{s \to s_j} |\bar{H}(\gamma_N(s))|=0.
$$
Now we set $p_n:=\gamma_N(s)$, so we obtain the result.
\enD

\begin{theorem}\label{the Gaussian curvature of X}
Let $(\bX,\bv,\bw):U \to \R_1^3\times\Delta_4$ be a lightcone framed surface and $\bF:U\to\R_1^3$ be the focal surface of $\bX$.
Suppose that the Gaussian curvature $K$ of $\bX$ is bounded at each lightlike point $p$ of $\bX$.
Then there exists a sequence $\{p_n\}_{n=1,2,3,...}$ in $U_+$(resp. $U_-$) converging to $p$ such that 
$$
\lim_{n \to \infty} K(p_n)=0 \ or \ \lim_{n \to \infty} K(p_n)=p_0\ne0,
$$
where $p_0$ is a constant.
\end{theorem}
\demo
Let $q_1,q_2\in V\subset U$ satisfy ${\tilde{\lambda}}(q_1)>0,{\tilde{\lambda}}(q_2)<0$ and $\gamma(s),s\in[0,2\pi]$ be a smooth curve in $U$ satisfied $\gamma(0)=q_1,\gamma(2\pi)=q_2$.
Then by Fourier series expansion of $\gamma(s)$, there exists a sufficiently large positive integer $N$, the real analytic curve defined by $\gamma_N(s)$ satisfies
\begin{eqnarray}\label{lambda-gamma-2}
\tilde{\lambda}(\gamma_N(0))>0,\tilde{\lambda}(\gamma_N(2\pi))<0.
\end{eqnarray}
We set 
$$
\check{\lambda}(s):=\tilde{\lambda}(\gamma_N(s)),s\in[0,2\pi]
$$
Since $\check{\lambda}(s)$ is a real analytic function defined on $[0,2\pi]$, the set of zeros of the function $\check{\lambda}(s)$ consists of a finite set of points
$$
0<s_1<s_2<\cdot\cdot\cdot<s_n<2\pi.
$$
By \eqref{lambda-gamma-2}, we can choose the number $j$ such that the sign of $\check{\lambda}(s)$ changes from positive to negative at $s=s_j$.
So there exists a positive integer $m$ such that
$$
\lim_{s \to s_j}\frac{\check{\lambda}(s)}{(s-s_j)^m}=\lambda_0\ne0,
$$
where $\lambda_0$ is a non-zero real number.
By Remark \ref{NLMNKH extend to U}, $\hat{K}$ can be extended to $U$.
We set
$$
\check{\hat{K}}(s):=\hat{K}(\gamma_N(s)),s\in[0,2\pi],
$$
then 
$$
K(\gamma_N(s)):=\frac{\check{\hat{K}}(s)}{|\check{\lambda}(s)|^2},s\ne s_1,s_2,\cdot\cdot\cdot,s_n.
$$
As $K(p)$ is bounded at each $p\in L(\bX)$, we get $\hat{K}(s_j)=0$.
Since $\hat{K}(s)=0$ is a real analytic function, there exists a positive integer $l$ such that 
$$
\lim_{s \to s_j} \frac{\check{\hat{K}}(s)}{(s-s_j)^l}=K_0\ne0,
$$
where $K_0$ is a non-zero real number.
Then it holds that
\begin{eqnarray*}
K(\gamma_N(s))=\frac{\check{\hat{K}}(s)}{(s-s_j)^l} \frac{(s-s_j)^{2m}}{|\check{\lambda}(s)|^{2m}}(s-s_j)^{l-2m},
\end{eqnarray*}
\begin{eqnarray*}
\lim_{s \to s_j} |s-s_j|^{2m-l}|K(\gamma_N(s))|=\frac{|K_0|}{|\lambda_0|^2}\ne0.
\end{eqnarray*}
As $K(p)$ is bounded, we have $l-2m\ge0$.
If $l-2m>0$, we have 
$$
\lim_{s \to s_j} |K(\gamma_N(s))|=0;
$$
 if $l-2m=0$, we have 
$$
\lim_{s \to s_j}|K(\gamma_N(s))|=\frac{|K_0|}{|\lambda_0|^2}\ne0.
$$
Now we set $p_n:=\gamma_N(s)$, so we obtain the result.
\enD

\begin{corollary}\label{the Gaussian curvature of F} 
Let $(\bX,\bv,\bw):U \to \R_1^3\times\Delta_4$ be a lightcone framed surface and $\bF:U\to\R_1^3$ be the focal surface of $\bX$.
Suppose that the Gaussian curvature $K$ of $\bX$ and $\bar{K}$ of $\bF$ are both bounded at each lightlike point of $\bX$.
Then there exists a sequence $\{p_n\}_{n=1,2,3,...}$ in $U_+$(resp. $U_-$) converging to $p$ such that 
$$
\lim_{n \to \infty} \bar{K}(p_n)=0 \ or \ \lim_{n \to \infty} \bar{K}(p_n)=\bar{p}_0\ne0,
$$ where $\bar{p}_0$ is a constant.
\end{corollary}
\demo
Let $q_1,q_2\in V\subset U$ satisfy $\tilde{\bar{\lambda}}(q_1)>0,\tilde{\bar{\lambda}}(q_2)<0$ and $\gamma(s),s\in[0,2\pi]$ be a smooth curve in $U$ satisfied $\gamma(0)=q_1,\gamma(2\pi)=q_2$.
Then by Fourier series expansion of $\gamma(s)$, there exists a sufficiently large positive integer $N$, the real analytic curve defined by $\gamma_N(s)$ satisfies
\begin{eqnarray}\label{lambda-gamma-3}
\tilde{\bar{\lambda}}(\gamma_N(0))>0,\tilde{\bar{\lambda}}(\gamma_N(2\pi))<0.
\end{eqnarray}
We set 
$$
\check{\bar{\lambda}}(s):=\tilde{\bar{\lambda}}(\gamma_N(s)),s\in[0,2\pi]
$$
Since $\check{\bar{\lambda}}(s)$ is a real analytic function defined on $[0,2\pi]$, the set of zeros of the function $\check{\bar{\lambda}}(s)$ consists of a finite set of points
$$
0<s_1<s_2<\cdot\cdot\cdot<s_n<2\pi.
$$
By \eqref{lambda-gamma-3}, we can choose the number $j$ such that the sign of $\check{\bar{\lambda}}(s)$ changes from positive to negative at $s=s_j$.
So there exists a positive integer $m$ such that
$$
\lim_{s \to s_j}\frac{\check{\bar{\lambda}}(s)}{(s-s_j)^m}=\bar{\lambda}_0\ne0,
$$
where $\bar{\lambda}_0$ is a non-zero real number.
By Remark \ref{NLMNKH extend to U}, $\hat{K}$ can be extended to $U$.
We set
$$
\check{\hat{\bar{K}}}(s):=\hat{\bar{K}}(\gamma_N(s)),s\in[0,2\pi],
$$
then 
$$
\bar{K}(\gamma_N(s)):=\frac{\check{\hat{\bar{K}}}(s)}{|\check{\bar{\lambda}}(s)|^2},s\ne s_1,s_2,\cdot\cdot\cdot,s_n.
$$
As $\bar{K}(p)$ is bounded at each $p\in L(\bX)$, we get $\hat{\bar{K}}(s_j)=0$.
Since $\hat{\bar{K}}(s)=0$ is a real analytic function, there exists a positive integer $l$ such that 
$$
\lim_{s \to s_j} \frac{\check{\hat{\bar{K}}}(s)}{(s-s_j)^l}=\bar{K}_0\ne0,
$$
where $\bar{K}_0$ is a non-zero real number.
Then it holds that
\begin{eqnarray*}
&&\bar{K}(\gamma_N(s))=\frac{\check{\hat{\bar{K}}}(s)}{(s-s_j)^l} \frac{(s-s_j)^{2m}}{|\check{\bar{\lambda}}(s)|^2}(s-s_j)^{l-2m},\\
&&\lim_{s \to s_j} |s-s_j|^{2m-l}|\bar{K}(\gamma_N(s))|=\frac{|\bar{K}_0|}{|\bar{\lambda}_0|^2}\ne0.
\end{eqnarray*}
As $\bar{K}(p)$ is bounded, we have $l-2m\ge0$.
If $l-2m>0$, we have 
$$
\lim_{s \to s_j} |\bar{K}(\gamma_N(s))|=0;
$$
 if $l-2m=0$, we have 
$$
\lim_{s \to s_j}|\bar{K}(\gamma_N(s))|=\frac{|\bar{K}_0|}{|\bar{\lambda}_0|^2}\ne0.
$$
Now we set $p_n:=\gamma_N(s)$, so we obtain the result.
\enD

\section{Example}

\begin{example}
We define $(\bX,\bv,\bw):[0,2\pi]\times[0,2\pi] \to \R^3_1 \times \Delta_4$ by
\begin{eqnarray*}
&&\bX=(\sin u, \cos u \sin v, \cos u \cos v),\\
&&\bv=(1, \sin v, \cos v),\\
&&\bw=(1,-\sin v,-\cos v).
\end{eqnarray*}
Then $(\bX,\bv,\bw):U \to \R^3_1 \times \Delta_4$ is a lightcone framed surface with lightcone frame. 
By a direct calculation, we have
\begin{eqnarray*}
\bm=(0,-\cos v,\sin v).
\end{eqnarray*}
and the curvature $(A,B,C)=(0,\frac{1}{2}\cos u(\sin u-cos u),-\frac{1}{2}\cos u(\sin u+cos u))$.
Then $\{\bv,\bw,\bm\}$ as above is the lightcone frame of $\bX$.
The basic invariants $(\mathcal G,\mathcal F_1,\mathcal F_2)$ of $\bX$ are as follows:
\begin{eqnarray*}
&&\mathcal G=\begin{pmatrix}
		a_1&  b_1& c_1\\
		a_2&  b_2& c_2
	\end{pmatrix}=\begin{pmatrix}
		-\frac{1}{2}(\sin u-\cos u)&  \frac{1}{2}(\sin u+\cos u)& 0\\
		0&  0& -\cos u
	\end{pmatrix},\\
&&\mathcal F_1=\begin{pmatrix}
		e_1&  0& 2g_1\\
		0& -e_1& 2f_1\\
		f_1& g_1& 0
	\end{pmatrix}=\begin{pmatrix}
		0&  0& 0\\
		0& 0& 0\\
		0& 0& 0
	\end{pmatrix},\\
&&\mathcal F_2=\begin{pmatrix}
		e_2&  0& 2g_2\\
		0& -e_2& 2f_2\\
		f_2& g_2& 0
	\end{pmatrix}=\begin{pmatrix}
		0&  0& -1\\
		0& 0& 1\\
		\frac{1}{2}& -\frac{1}{2}& 0
	\end{pmatrix}.
\end{eqnarray*}

The lightlike point set and singular point set of $\bX$ are 
\begin{eqnarray*}
&&L(\bX)=\{ (u,v)\in [0,2\pi]\times[0,2\pi] \ | \ u=k\pi+\frac{\pi}{4},k=0,1 \},\\
&&S(\bX)=\{ (u,v)\in [0,2\pi]\times[0,2\pi] \ | \ u=k\pi+\frac{\pi}{2},k=0,1 \}.
\end{eqnarray*}

Next we give the case for $c_2>0$, that is $u\in (\frac{\pi}{2},\frac{3\pi}{2})$.
The case for $c_2<0$ is simiar to the case for $c_2>0$.
By a direct calculation, we have:
\begin{eqnarray*}
&&\tilde{\lambda}=-\cos 2u,\\
&&\hat{\bn}=(\sin u,-\cos u\sin v,-\cos u\cos v),\\
&&\hat{K}=-\cos u,\\
&&\hat{H}=-\sin^2u\cos u,\\
&&\mu=-1 \ or \ \cos 2u.
\end{eqnarray*}
Then the focal surface of the lightcone framed surface $\bX$ and invariants of $\bF$ are as follows: \\
$(1)$ $\mu=-1$, 
\begin{eqnarray*}
&&\bF=(2\sin u,0,0),\\
&&\bar{a}_1=-\frac{1}{2}(\sin u-\cos u),\\
&&\bar{b}_1=\frac{1}{2}(\sin u+\cos u),\\
&&\bar{c}_2=-\cos u,\\
&&\bar{a}_2=\bar{b}_2=\bar{c}_1=0.
\end{eqnarray*}
$(2)$ $\mu=\cos 2u$, 
\begin{eqnarray*}
&&\bF=(2\sin^3 u,2\cos^3 u\sin v,2\cos^3 u\cos v),\\
&&\bar{a}_1=-\frac{1}{2}\sin u(2\cos^2 u+\sin 2u)+\frac{1}{2}\cos u(2\sin^2 u-\sin 2u),\\
&&\bar{b}_1=\frac{1}{2}\sin u(2\cos^2 u+\sin 2u)+\frac{1}{2}\cos u(2\sin^2 u+\sin 2u),\\
&&\bar{c}_2=-2\cos^3 u,\\
&&\bar{a}_2=\bar{b}_2=\bar{c}_1=0.
\end{eqnarray*}
It is obvious that $\bF$ is a lightcone framed surface with the lightcone frame $\{\bv,\bw,\bm\}$.

\begin{center}
 \centering
 \begin{minipage}[c]{0.4\textwidth}
  \centering
  \includegraphics[scale=0.3]{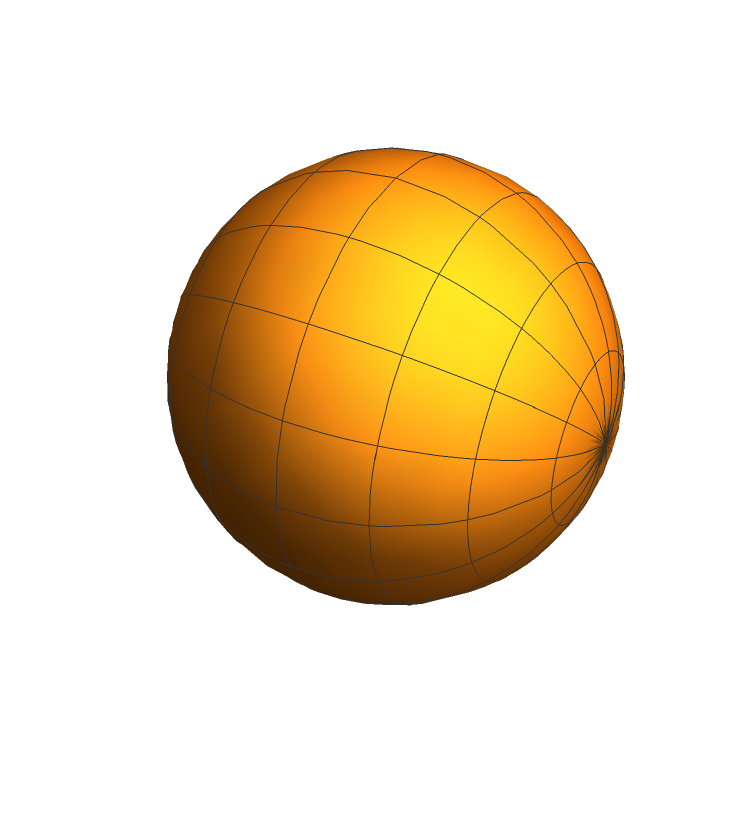} \\
   
 \end{minipage}
\vspace*{2mm}\\ Figure 1\quad the lightcone framed surface
\end{center}

\begin{center}
 \centering
 \begin{minipage}[c]{0.4\textwidth}
  \centering
  \includegraphics[scale=0.8]{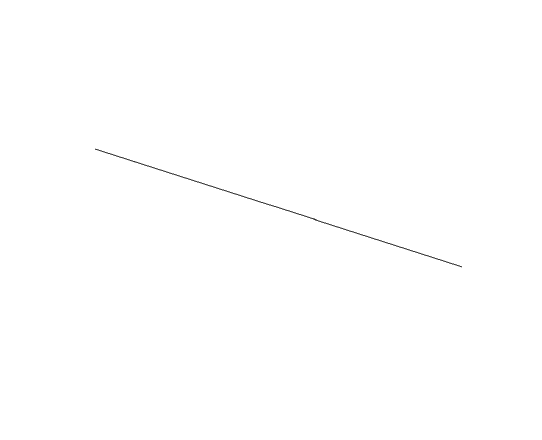} \\
   (i) $\mu=-1$
 \end{minipage}%
 \begin{minipage}[c]{0.4\textwidth}
  \centering
  \includegraphics[scale=0.5]{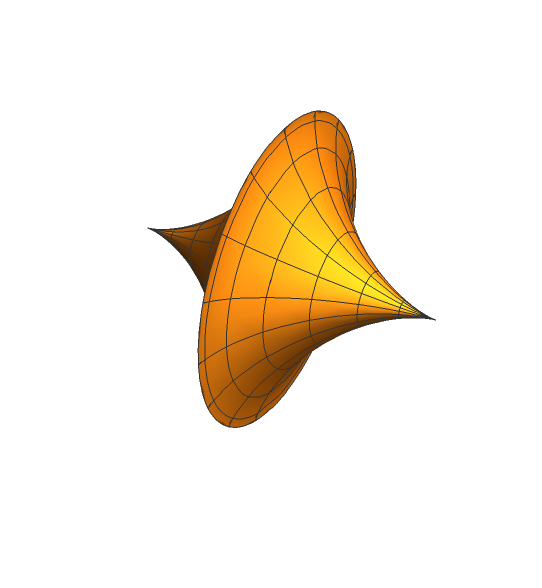} \\
 (ii) $\mu=\cos 2u$\\
 \end{minipage}
\vspace*{2mm}\\ Figure 2\quad the focal surface of the lightcone framed surface
\end{center}

\end{example}


\ \\
Chang Xu,
\\
School of Mathematics and Statistics, Northeast Normal University, Changchun 130024, P.R.CHINA
\\
E-mail address: xuc052@nenu.edu.cn
\\
\\
Liang Chen,
\\
School of Mathematics and Statistics, Northeast Normal University, Changchun 130024, P.R.CHINA
\\
E-mail address: chenl234@nenu.edu.cn


\begin{thebibliography}{99}
{\small


\bibitem{Arnold2} V.I.Arnol'd, S.M.Gusein-Zade and A.N.Varchenko,
\newblock{\em Singularities of Differentiable Maps vol. {\rm I}}.
\newblock Birkh\"auser (1986). 


\bibitem{Bruce-Giblin} J.W.Bruce and P.J.Giblin, 
\newblock{\em Curves and singularities. A geometrical introduction to singularity theory. Second edition}. 
\newblock{Cambridge University Press}, Cambridge (1992).


\bibitem{Liang Chen1} L.Chen and M.Talahashi, 
\newblock{Lightcone framed curves in the Lorentz-Minkowski $3$-space}. 
\newblock{\em Turkish J.Math}, {\bf 48} (2024), 307-326.


\bibitem{T.Liu1} T.Liu and D.Pei,
\newblock{Mixed type curves and the lightcone frame in Minkowski $3$-space.} 
\newblock{\em International Journal of Geometric Methods in Modern Physics}, {\bf 17} (2020) 2050088, 14.


\bibitem{M.Li1} M.Li, D.Pei and M.Takahashi,
\newblock{Lightcone framed surfaces in Lorentz-Minkowski $3$-space.} 
\newblock{\em preprint}, 2024.


\bibitem{Tari1} F.Tari,
\newblock{Caustics of surfaces in the Minkowski 3-space}. 
\newblock{\em Quaet.J.Math}, {\bf 63} (2012), 189-209.


\bibitem{Fukunaga1} T.Fukunaga and M.Takahashi, 
\newblock{Framed surfaces in Euclidean Space}. 
\newblock{\em Bull Braz Math Soc, New Series}, {\bf 50} (2019), 37-65.


\bibitem{Honda1} S.Honda and M.Takahashi,
\newblock{Evolutes and focal surfaces of framed immersions in the Euclidean space.} 
\newblock{\em Proceedings of the Royal Society of Edinburgh}, {\bf 150} (2020) 497--512.


\bibitem{Honda1} A.Honda, M.Koiso, M.Kokubu, M.Umehara and K.Yamada, 
\newblock{Mixed type surfaces with bounded mean curvature in $3$-dimensional space-times}. 
\newblock{\em Differential Grometry and its Applications}, {\bf 52} (2017), 64-77.


\bibitem{Honda2} A.Honda, K.Saji and K.Teramoto, 
\newblock{Mixed type surfaces with bounded Gaussian curvature in three-dimensional Lorentzian manifolds}. 
\newblock{\em Advance in Mathematics}, {\bf 365} (2020), 107036.


\bibitem{Teramoto1} K.Teramoto,
\newblock{Principal curvatures and parallel surfaces of wave fronts}. 
\newblock{\em Adv,Geom}, {\bf 19} (2019), 541-554.


\bibitem{Teramoto2} K.Teramoto, 
\newblock{Focal surfaces of wave fronts in the Euclidean $3$-space}. 
\newblock{\em Glasgow Math.J}, {\bf 61} (2019), 425-440.


\bibitem{Teramoto3} K.Teramoto, 
\newblock{Focal surfaces of fronts associated to unbounded principal curvatures}. 
\newblock{\em Rocky Mountain Journal of Mathematics}, vol.53, {\bf 5} (2023), 1587-1608.


\bibitem{Teramoto4} K.Teramoto,
\newblock{Parallel and focal surfaces of wave fronts}. 
\newblock{\em Kobe university repository}, (2023).


\bibitem{Teramoto5} K.Teramoto,
\newblock{Parallel and dual surfaces of cuspidal edges}. 
\newblock{\em Differ.Geom.Appl}, {\bf 44} (2016), 52-62.


\bibitem{Bruce1} J.W.Bruce and P.J.Giblin, 
\newblock{Families of surfaces: focal sets, ridges and umbilics}. 
\newblock{\em Math. Oroc. Camb. Soc}, {\bf 125} (1999), 243.


\bibitem{Fukui1} T.Fukui and M.Hasegawa, 
\newblock{Singularities of parallel surfaces}. 
\newblock{\em Tohoku Math.J.}, {\bf 64} (2012), 387-408.


\bibitem{Martins1} L.F.Martins, K.Saji, M.Umehara and K.Yamada, 
\newblock{Behavior of Gaussian curvature and Mean curvature near non-degenerate singular points on wave fronts}. 
\newblock{\em Geometry and topology of manifolds}, vol.154, (2016), 445-462.


\bibitem{Saji1} K.Saji, M.Umehara and K.Yamada, 
\newblock{Behavior of corank one singular points on wave fronts}. 
\newblock{\em Kyushu J.Math}, {\bf 62} (2008), 259-280.


\bibitem{Saji2} K.Saji, 
\newblock{Criteria for $D_4$ singularities of wave fronts}.
\newblock{\em Tohoku Math.J.}, {\bf 63} (2011), 137-147.


\bibitem{Fujimori1} S.Fujimori, K.Saji and M.Umehara, 
\newblock{Singularities of maximal surfaces}. 
\newblock{\em Math.Z.}, {\bf 259} (2008), 827-848.


\bibitem{Izumiya1} S.Izumiya and K.Saji, 
\newblock{The mandala of Legendrian dualities for pseudo-spheres in Lorentz-Minkowski space and "flat" spacelike surfaces}. 
\newblock{\em Journal of Singularities}, vol.2, (2010), 92-127.


}
\end{thebibliography}
\end{document}